\documentclass[a4paper]{amsart}
\usepackage{srcltx}
\usepackage{amsmath}
\usepackage{amsthm}
\usepackage{amsfonts}
\usepackage{amssymb}

\usepackage{xr}
\newcommand\ECref[1]{\cite[\ref{enlcat-#1}]{enlcat}}
\newcommand\ECrefs[2]{\cite[\ref{enlcat-#1}\ref{enlcat-#2}]{enlcat}}

\newlength{\abstand}
\def\F{{\mathbb F}}
\def\G{{\mathbb G}}
\def\R{{\mathbb R}}
\def\C{{\mathbb C}}
\def\N{{\mathbb N}}
\def\Nn{{{\mathbb N}_0}}
\def\Np{{{\mathbb N}_+}}
\def\Q{{\mathbb Q}}
\def\Z{{\mathbb Z}}
\def\P{{\mathbb P}}
\def\G{{\mathbb G}}
\def\Ns{\str{\mathbb N}}
\def\Nns{\str{{\mathbb N}_0}}
\def\Nps{\str{{\mathbb N}_+}}
\def\Cs{\str{\mathbb C}}
\def\Ps{\str{\mathbb P}}
\def\Qs{\str{\mathbb Q}}
\def\Rs{\str{\mathbb R}}
\def\Zs{\str{\mathbb Z}}
\def\Qab{{\Q^{\text{ab}}}}
\def\Qb{\bar{\Q}}
\def\ordp{{\text{ord}_p}}
\def\p{{\mathfrak p}}
\def\pinf{{(p^\infty)}}
\def\resp{resp.\ }
\def\Ob{\mathrm{Ob}}
\def\cA{\mathcal A}
\def\cB{\mathcal B}
\def\cC{\mathcal C}
\def\cD{\mathcal D}
\def\cI{\mathcal I}
\def\cJ{\mathcal J}
\def\cL{\mathcal L}
\def\cP{\mathcal P}
\def\cR{\mathcal R}
\def\cF{\mathcal F}
\def\cG{\mathcal G}
\def\cO{\mathcal O}
\def\RR{{\bf{R }}}

\catcode`\ä=\active
\def ä{\"a}
\catcode`\ö=\active
\def ö{\"o}
\catcode`\ü=\active
\def ü{\"u}
\catcode`\Ä=\active
\def Ä{\"A}
\catcode`\Ö=\active
\def Ö{\"O}
\catcode`\Ü=\active
\def Ü{\"U}
\catcode`\ß=\active
\def ß{\ss{}}
\catcode`\é=\active
\def é{\'{e}}
\catcode`\É=\active
\def É{\'{E}}

\def\ra{\rightarrow}

\newcommand\str[1]{{\mbox{}^*#1}}
\newcommand\Mor[3]{\mathrm{Mor}_{#1}(#2,#3)}
\newcommand\Hom[3]{\mathrm{Hom}_{#1}(#2,#3)}

\newcommand\cov[2]{\mathrm{Funct}\,(#1,#2)}

\def\Ensfin{{\mathrm{Ens}}^{\mathrm{fin}}}

\newcommand\fin[1]{#1^{\text{fin}}}
\def\finG{\fin{\tilde{G}}}
\def\strfin{$\str{\text{finite}}$}

\swapnumbers
\theoremstyle{definition}
\newtheorem{defi}{Definition}[section]
\newtheorem{bsp}[defi]{Example}
\newtheorem{bspe}[defi]{Examples}
\newtheorem{satzdefi}[defi]{Proposition/ Definition}
\newtheorem{lemma}[defi]{Lemma}

\newtheorem{bem}[defi]{Remark}
\newtheorem{satz}[defi]{Proposition}
\newtheorem{thm}[defi]{Theorem}
 \newtheorem{cor}[defi]{Corollary}

\input xy
\xyoption{all}

\title{Nonstandard étale cohomology}
\author{Lars Brünjes, Christian Serp\'{e}}
\date{\today}
\thanks{The first author is supported by the European Union.}
\thanks{This work has been written under hospitality of the Department of Pure
  Mathematics and Mathematical Statistics of the University of Cambridge and
  the "Sonderforschungsbereich geometrische Strukturen in der Mathematik"
  of the University of Münster}
\address{
  Lars Brünjes \\
  Universität Regensburg \\
  NWF I - Mathematik\\
  93040 Regensburg \\
  Germany}
\email{Lars.Bruenjes@mathematik.uni-regensburg-de}
\address{
  Christian Serp\'{e} \\
  Westfälische Wilhelms-Universität Münster \\
  Mathematisches Institut \\
  Sonderforschungbereich 478 ``Geometrische Strukturen in der Mathematik'' \\
  Hittorfstr. 27 \\
  D-48149 Münster \\
  Germany}
\email{serpe@math.uni-muenster.de}


\date{\today}

\begin{document}

\begin{abstract}
  A lot of good properties of étale cohomology only hold for torsion coefficients.
We use ``enlargement of categories" as developed in \cite{enlcat}
to define a cohomology theory that inherits the important properties of étale cohomology
while allowing greater flexibility with the coefficients.
In particular,
choosing coefficients $\str{\Z}/P\str{\Z}$
(for $P$ an infinite prime and $\str{\Z}$ the enlargement of $\Z$)
gives a Weil cohomology,
and choosing $\str{\Z}/l^h\str{\Z}$ (for $l$ a finite prime and $h$ an infinite number)
allows comparison with ordinary $l$-adic cohomology.
More generally,
for every $N\in\str{\Z}$,
we get a category of $\str{\Z}/N\str{\Z}$--constructible sheaves
with good properties.

\end{abstract}

\maketitle
\tableofcontents

\vspace{\abstand}

\section{Introduction}

\mbox{}\\
\vspace{\abstand}

Étale cohomology is among the most important tools
of modern algebraic and arithmetic geometry,
and one of the main reasons for its importance
is the fact
that étale cohomology
allows us to construct a ``good" cohomology theory with coefficients in a field of characteristic zero
for varieties over fields of arbitrary characteristic --- the $l$-adic cohomology.

But even though étale cohomology
is defined for coefficients in arbitrary étale sheaves,
we only get good properties (like proper or smooth base change, finiteness theorems or Poincaré duality),
if we restrict ourselves to torsion sheaves;
often it is even necessary to only consider constructible sheaves
or --- more restrictive still --- (locally) constant constructible sheaves.

Therefore, in order to define $l$-adic cohomology of a variety $X$,
we can not simply take étale cohomology of $X$ with coefficients in the constant sheaf $\Q_l$ (which is not torsion),
but have to take $\varprojlim H^i(X,\Z/l^n)\otimes_{\Z_l}\Q_l$ instead.

As a consequence,
$l$-adic cohomology is not in general a derived functor,
and not all properties of étale cohomology carry over to $l$-adic cohomology.

For example, if $X$ is a variety over a field $k$,
if $k_s$ is a separable closure of $k$,
and if $\cF$ is an étale sheaf on $X$,
then we have the Hochschild-Serre spectral sequence
$H^p(Gal(k_s/k),H^q(X\otimes_kk_s,\cF))\Rightarrow H^{p+q}(X,\cF)$,
but we do not in general have such a spectral sequence when we consider $l$-adic cohomology instead.

Using an ``honest" derived functor can overcome these disadvantages,
as was for example demonstrated by Jannsen who introduced \emph{continuous étale cohomology} in \cite{uwe}
and showed that it coincides with $l$-adic cohomology for projective varieties
over separably closed fields.\\

In this paper,
we use the theory of ``enlargement of categories" as developed in \cite{enlcat}
to define \emph{*étale cohomology},
as the enlargement of ordinary étale cohomology
(by results from \cite{enlcat}, *étale cohomology is then in particular a derived functor);
this allows us to choose coefficients in arbitrary *étale sheaves.

By the transfer principle,
the good properties of étale cohomology
for torsion sheaves,
constructible and (locally) constant constructible sheaves
carry over to analogous properties of *étale cohomology
for *torsion sheaves,
*constructible and (locally) constant *constructible sheaves.

An example of a constant *constructible sheaf is the *sheaf associated to
$\F_P:=\str{\Z}/P\str{\Z}$ for an infinite prime $P$.
Internally, $\F_P$ is a *finite field,
but externally, $\F_P$ is a field of characteristic zero.
Therefore *étale cohomology with coefficients in $\F_P$ is a cohomology theory
with coefficients in a field of characteristic zero (it is in fact a Weil cohomology)
which has all the good properties of étale cohomology with torsion coefficients ---
including, for example, a Hochschild-Serre spectral sequence
$H^p(\str{Gal(k_s/k)},H^q(X\otimes_kk_s,\F_P))\Rightarrow H^{p+q}(X,\F_P)$ in the situation described above.

This cohomology with coefficients in $\F_P$
is closely related to the one constructed by Toma\v{s}i\'{c} (\cite{tomasic}) using ultrafilters
(for a similar construction also compare Serre \cite{serre2}),
and the proof that it is indeed a Weil cohomology is the same as Toma\v{s}i\'{c}'s proof.
But while Toma\v{s}i\'{c}'s approach is more direct,
ours has the advantage of giving greater flexibility in the choice of coefficients
and of producing corresponding categories of *étale sheaves as well.

The *sheaf associated to
$\str{\Z}/l^h\str{\Z}$ for an ordinary prime $l$ and an infinite natural number $h$
is another example of a constant *constructible sheaf.
Here the natural projections $\str{\Z}/l^h\str{\Z}\rightarrow\str{\Z}/l^n\str{\Z}\cong\Z/l^n\Z$
often allow us to compare
*étale cohomology with coefficients in $\str{\Z}/l^h\str{\Z}$
to ordinary $l$-adic cohomology,
and we get a close correspondence between *constructible $\str{\Z}/l^h\str{\Z}$--sheaves
and $l$-adic sheaves.

So while inheriting all the good properties of étale cohomology,
we gain much greater flexibility concerning suitable coefficients
when using *étale cohomology.

In addition to that,
*étale cohomology is not only defined for ordinary schemes,
but for *schemes as well,
and even though we will not explore that alley in the present paper,
it seems to open a lot of promising possibilities.\\



The plan of the paper is as follows:
In the second chapter, we define the relevant categories of schemes and étale sheaves
and their enlargements,
and we define *étale cohomology and *étale cohomology with compact support
as the enlargement of ordinary étale cohomology respectively étale cohomology with compact support.
We show that the usual properties hold in the nonstandard setup as well;
in particular we get the familiar Leray spectral sequence.

In the third chapter, we specialize to the case of varieties over a field. We define geometric and arithmetic
*étale cohomology, we show that we have a Hochschild-Serre spectral sequence as mentioned above,
prove Poincaré duality for *étale cohomology
as well as the existence of a cycle map,
and then formulate a meta-theorem about the transfer of properties from ordinary étale cohomology
to *étale cohomology.
Further we prove that
if we take *étale cohomology with coefficients in the constant *sheaf $\str{\Z}/P\str{\Z}$,
the resulting cohomology groups are \emph{finite} (and not just *finite),
and we see that we get a Weil cohomology in that way.

In the fourth chapter, we specialize even further to the case of smooth projective varieties over $\C$.
We show that *étale cohomology with coefficients in a *finite constant sheaf $M$
is then nothing else than
ordinary singular cohomology with coefficients in $M$.

In the fifth chapter, we investigate the connection between $l$-adic sheaves and
*constructible $\str{\Z}/l^h$--sheaves,
and we prove a comparison theorem between geometric *étale cohomology
and geometric $l$-adic cohomology.

In the appendix, we recall the most important definitions from the theory of
enlargements and enlargements of categories.\\

We thank Niko Naumann and Christopher Deninger for helpful discussions.

\vspace{\abstand}

\newcommand\spec[1]{{\text{Spec}\,\left(#1\right)}}
\newcommand\etshv[1]{{\text{Shv}_{et}^{Ab}\left(#1\right)}}
\newcommand\preshv[1]{{\text{PreShv}^{Ab}\left(#1\right)}}
\newcommand\etsh{{\text{Shv}_{et}^{Ab}}}
\newcommand\constr[1]{{\text{Constr}_{et}\left(#1\right)}}
\newcommand\locconst[1]{{\text{LocConst}_{et}\left(#1\right)}}
\newcommand\modconstr[2]{{\text{ModConstr}_{et}^{Ab}\left(#1,#2\right)}}
\newcommand\Mod[1]{{\text{Mod}\left(#1\right)}}
\newcommand\modshv[2]{{\text{ModShv}_{et}^{Ab}\left(#1,#2\right)}}
\newcommand\modshvs[1]{{\text{ModShv}_{et}^{Ab}\left(#1\right)}}

\newcommand\derconstr[2]{\cD^b_{constr}\left(#1,#2\right)}
\newcommand\derconstrs[1]{\cD^b_{constr}\left(#1\right)}
\newcommand\derconstrsf[1]{\cD^{fb}_{constr}\left(#1\right)}
\newcommand\derconstrf[2]{\cD^{fb}_{constr}\left(#1,#2\right)}
\def\finring{\mathrm{FinRings}}
\def\Abe{\mathrm{Ab}}
\def\finab{\fin{\Abe}}

\def\etshvo{{\text{Shv}_{et}^{Ab}}}
\def\constro{{\text{Constr}_{et}}}
\def\locconsto{{\text{LocConst}_{et}}}
\def\fib{\mathcal Fib}

\section{Enlargement of schemes and étale sheaves}\label{chone}

\vspace{\abstand}

In the first part of this chapter we introduce the relevant categories
of schemes and their enlargements. For the applications we have in mind
in this paper,
it is not strictly necessary to enlarge the category of
schemes as well,
but for us it seems most natural to view the situation in this full generality.
In the second part
we consider two classes of fibrations:
One is étale sheaves over schemes,
the other one is étale sheaves over finite rings
(which will be most relevant for the following chapters).
The study of their enlargements will enable us to define nonstandard cohomology and prove its elementary properties
in the third part of this chapter.
Finally, in the last part, we define nonstandard cohomology with compact support.

Most parts of this chapter are an easy
application of \cite{enlcat} where the authors developed  a general theory of
enlargement of categories.

\mbox{}\\
\vspace{\abstand}

\subsection{Enlargement of schemes}

\mbox{}\\
\vspace{\abstand}

Let $S$ be a noetherian scheme. By $Sch/S$ we denote the category of schemes over $S$ which
are separated and of finite type over $S$. The category $Sch/S$ is small.
We choose a superstructure $\hat{S}$ such that the category $Sch/S$ is  $\hat{S}$--small.
Let $*:\hat{S} \ra\widehat{\str{S}}$ be
an enlargement. Now we want to apply the enlargement of categories which has been constructed
and whose properties have been studied in \cite{enlcat}.

\vspace{\abstand}

We call the objects of the category $\str{Sch/S}$  \emph{nonstandard schemes}  or \emph{*schemes over $\str{S}$}.
The objects in this category are not schemes, but can be thought of as
some kind of limit of schemes. The category $\str{Sch/S}$ has the
same formal properties as $Sch/S$.
For example, by \ECref{corfinlims} the category $\str{Sch/S}$ has a final
object and fibred products. We have the canonical faithful functor
$$*: Sch/S \ra \str{Sch/S},$$
and we call the objects in the image of this functor \emph{standard schemes}.

\vspace{\abstand}

Let $P$ be a property that a morphism of schemes can have. We assume that the identities are in $P$
and that $P$ is closed under composition. Then by $(Sch/S)^P$ we denote the subcategory of $Sch/S$
which has the same objects as $Sch/S$ but only those morphisms that have the property $P$. Then
$\str{(Sch/S)^P}$ is a subcategory of $\str{(Sch/S)}$, and we say that a morphism in $\str{(Sch/S)}$
has the property $P$ if the morphism lies in the subcategory $\str{(Sch/S)^P}$. We have the commutative
diagram

\vspace{\abstand}
$$
\xymatrix{
  (Sch/S)^P \ar@{^{(}->}[r] \ar@{^{(}->}[d] & (Sch/S) \ar@{^{(}->}[d] \\
  \str{(Sch/S)^P}  \ar@{^{(}->}[r]           & \str{(Sch/S)}
}
$$

\vspace{\abstand}

\noindent
and a morphism $f$ in $Sch/S$ has property $P$ iff $\str{f}$  has
property $P$ in $\str{(Sch/S)}$.

\vspace{\abstand}

\begin{bsp}
  In $Sch/S$ we have the properties of being an open immersion, of being proper
  and of being a compactifiable morphism. Here we say that a morphism $f:X\ra Y$
  in $(Sch/S)$ is compactifiable if there is another scheme $\bar X\in(Sch/S)$,
  an open immersion $i:X\hookrightarrow \bar X$, and a proper morphism
  $\bar f: \bar X\ra Y$ such that $f=\bar f \circ i$. By the transfer principle
  it then follows that a morphism $f:X\ra Y$ in $\str{(Sch/S)}$ is compactifiable
  if and only if there is a scheme $\bar X \in \str{(Sch/S)}$, an open immersion
  $i:X \hookrightarrow \bar X$ in $\str{(Sch/S)}$, and a proper morphism $p:\bar X\ra X$ in
  $\str{(Sch/S)}$ such that $f=\bar f\circ i$ holds. We denote the category of
  compactifiable morphisms by $(Sch/S)^{compac.}$.
\end{bsp}

\vspace{\abstand}

\subsection{Enlargement of étale sheaves}

\mbox{}\\
\vspace{\abstand}

For a scheme $X\in Sch/S$ we denote the small étale site of $X$ by $Et(X)$.
We choose once and for all an universe $U$ in such a way that for all $X\in Sch/S$ the category
$$\text{PreShv}(Et(X)):=\text{PreShv}(Et(X), U-Sets):=\cov{(Et(X))^{op}}{U-Sets}$$
contains all representable presheaves. For all $X \in Sch/S$ we denote by
$$\preshv{Et(X)}:=\text{PreShv}(Et(X), U-Ab):=\cov{(Et(X))^{op}}{U-Ab}$$
the category of abelian group objects in $\text{PreShv(Et(X))}$. This is
again a small category. We further denote by
$$\etshv{X}\subset \preshv{Et(X)}$$
the full subcategory of abelian sheaves in the étale topology.
We denote by $\fib_{\etsh}$ the following category. The objects are
pairs $(X,\cF)$ with $X\in Sch/S$ and $\cF\in\etshv{X}$.
A morphism from $(X,\cF)$ to $(Y,\cG)$ in $\fib_{\etsh}$ is a pair $(f,\varphi)$, where
$f:X\ra Y$ is a morphism of schemes in $(Sch/S)$, and $\varphi:f^*\cG\ra\cF$ is
a morphism in $\etshv{X}$. The obvious functor
$$ p: \fib_{\etsh} \ra (Sch/S)$$
which maps a pair $(X,\cF)$ to $X$
and a morphism $(f,\varphi):(X,\cF)\ra(Y,\cG)$ to $f:X\ra Y$ is an abelian fibration.

\vspace{\abstand}

\begin{bem}
  The fiber $\fib_{\etsh}(X)$ of $p$ in $X$ is just $\etshv{X}^{op}$.
\end{bem}

\vspace{\abstand}

For all objects $X,Y\in Sch/S$ and all morphisms $f:X\ra Y$ there is
a right adjoint functor $f_*:\etshv{X}\ra\etshv{Y}$ of $f^*$ and therefore
$p:\fib_{\etsh}\ra (Sch/S)$ is also a cofibration.
The category $\fib_{\etsh}$ is again small. All categories of sheaves which we will consider
are more or less subcategories of $\fib_{\etsh}$ and therefore also small,
and we will not comment on this anymore.
For a scheme $X\in Sch/S$, we denote by $\constr{X}$ the
category of constructible sheaves on $X$ (an abelian
étale sheaf $\cF$ is in $\constr{X}$ if and only if $X$ can be written
as a finite union of locally closed subschemes $Y\subseteq X$
such that $\cF_{|Y}$ is finite and locally constant).
As above we get an abelian fibration of categories
$$ \fib_{\text{Constr}_{et}}\ra (Sch/S).$$
In the same way we get the fibration of categories
$$ \fib_{\text{LocConst}_{et}} \ra (Sch/S)$$
of locally constant constructible étale sheaves.
We get the following full subcategories of fibrations:

$$\xymatrix{
  \fib_{\text{LocConst}_{et}} \ar@{^{(}->}[r] \ar[rd] & \fib_{\text{Constr}_{et}} \ar@{^{(}->}[r] \ar[d]  &    \fib_{\etsh} \ar[ld] \\
                                             & Sch/S             }$$

\vspace{\abstand}

Now we want to introduce the second family of fibrations.
For that we denote by $\finring$ the category of finite commutative rings.
For each $X\in Sch/S$ and for each finite Ring $R\in\finring$ we have the category
$\modconstr{X}{R}$ of constructible étale sheaves of $R$--modules.
This again gives an abelian fibration
$$\text{ModConstr}_{et}(X) \ra (\finring)^{op}.$$
Objects in $\text{ModConstr}_{et}(X)$ are pairs $(\cF,R)$ where $\cF$ is
a constructible étale sheaf of $R$--modules on $X$. A morphism of two such pairs
is a pair $(f,\varphi):(\cF,R)\ra(\cG,S)$, where $\varphi:S\ra R$ is a morphism
of commutative rings and $f:\cG\otimes_S R\ra \cF$ is an $R$--morphism of étale sheaves.
Because the functor $\_\otimes_R S$ has a right adjoint,
the functor $\text{ModConstr}_{et}(X)\ra(\finring)^{op}$ is again a bifibration.
In the same way we get an abelian bifibration
$$\text{ModShv}_{et}(X)\ra(\finring)^{op}$$
where we restrict ourselves not only to constructible sheaves.

\vspace{\abstand}

\begin{bem}
  More generally one can show that there is a bifibration
  $$ \text{ModConstr}_{et}\ra (Sch/S)\times (\finring)^{op}$$
  such that the fibre of $(X,R)$ for an $S$--scheme $X$ and a finite
  commutative ring $R$ is just the category $(\modconstr{X}{R})^{op}$.
\end{bem}

\vspace{\abstand}

Now we want to enlarge these abelian fibrations. For that we assume
that our superstructure has been chosen in such a way
that all occurring categories are $\hat{S}$--small. By \ECref{abfib} the enlargement of
an abelian fibration is again an abelian fibration. So we get
the abelian fibrations:

$$\str{\fib_{\etsh}} \ra \str{(Sch/S)}$$
$$\str{\fib_{\text{LocConst}_{et}}} \ra \str{(Sch/S)}$$
$$\str{\fib_{\text{Constr}_{et}}}\ra \str{(Sch/S)}$$
and

$$\xymatrix{
  \str{\fib_{\text{LocConst}_{et}}} \ar@{^{(}->}[r] \ar[rd] & \str{\fib_{\text{Constr}_{et}}} \ar@{^{(}->}[r] \ar[d]  &    \str{\fib_{\etsh}} \ar[ld] \\
                                             & \str{Sch/S}             }$$

\noindent
For all $X\in\str{(Sch/S)}$ we set
$$\str{\etshv{X}}:=(\str{\fib_{\etsh}}(X))^{op}$$
$$\str{\locconst{X}}:=(\str{\fib_{\text{LocConst}_{et}}}(X))^{op}$$
$$\str{\constr{X}}:=(\str{\fib_{\text{Constr}_{et}}}(X))^{op}.$$
It follows further that $\str{p}:\str{\fib_{\etsh}}\ra\str{(Sch/S)}$ is also a cofibration.
In particular this means that for all $X,Y\in\str{(Sch/S)}$ and all morphisms $f\in\Hom{\str{(Sch/S)}}{X}{Y}$
we have pull backs
$$ f^*:\str{\etshv{Y}}\ra\str{\etshv{X}},$$
$$ f^*:\str{\locconst{Y}}\ra\str{\locconst{X}},$$
$$ f^*:\str{\constr{Y}}\ra\str{\constr{X}},$$
and a push forward
$$ f_*:\str{\etshv{X}}\ra\str{\etshv{Y}}$$
which is right adjoint to $f^*$.
We get

\vspace{\abstand}

\begin{satz}\label{elestar}
  For each nonstandard scheme $X\in\str{Sch/S}$, the three categories
  $\str{\locconst{X}}$,
  $\str{\constr{X}}$
  and
  $\str{\etshv{X}}$
  are abelian categories, and for each $X\in Sch/S$ we have
  \begin{itemize}
  \item $\str{\etshv{\str{X}}}=\str{(\etshv{X})}$
  \item $\str{\locconst{\str{X}}}=\str{(\locconst{X})}$
  \item $\str{\constr{\str{X}}}=\str{(\constr{X})}$.
  \end{itemize}
  Furthermore, for each $X\in\str{Sch/S}$, the abelian
  category $\str{\etshv{X}}$ has enough injective objects.
\end{satz}

\vspace{\abstand}

\begin{proof}
  This follows from \ECref{enlfib}  and \ECref{injfib}.
\end{proof}

\vspace{\abstand}

\begin{satz}
  We have canonical functors
  $$ *:\fib_{\etshvo}\ra \str{\fib_{\etshvo}}$$
  $$ *:\fib_{\locconsto}\ra\str{\fib_{\locconsto}}$$
  $$ *:\fib_{\constro}\ra\str{\fib_{\constro}}$$
  such that the following diagrams are commutative:
  $$\xymatrix{
    \fib_{\etshvo} \ar[r]^{*} \ar[d] & \str{\fib_{\etshvo}}\ar[d]  \\
    Sch/S \ar[r]^{*}          & \str{Sch/S}}$$
  $$\xymatrix{
    \fib_{\locconsto} \ar[r]^{*} \ar[d] & \str{\fib_{\locconsto}}\ar[d]  \\
    Sch/S \ar[r]^{*}          & \str{Sch/S}}$$
  $$\xymatrix{
    \fib_{\constro} \ar[r]^{*} \ar[d] & \str{\fib_{\constro}}\ar[d]  \\
    Sch/S \ar[r]^{*}          & \str{Sch/S}}$$
  If for an $X\in Sch/S$ we identify the category $\str{\etshv{\str{X}}}$
  with $\str{(\etshv{X})}$ by means of proposition \ref{elestar}, then the above functors
  become the usual functors $*:\etshv{X}\ra\str{(\etshv{X})}$, $*:\locconst{X}\ra\str{(\locconst{X})}$
  and $*:\constr{X}\ra\str{(\constr{X})}$.
\end{satz}

\vspace{\abstand}

\begin{proof}
  \ECrefs{satzfunctors}{sfunctorsiii} and \cite[\S7]{enlcat}.
\end{proof}

\vspace{\abstand}

\begin{bem}
  For $X\in Sch/S$, the objects of $\str{\etshv{X}}$ can be
  seen as so called internal functors from the category
  $\str{Et(X)}$ to the category $\str{(U-Ab)}$ which fulfil
  a certain sheaf-like condition --- for details we refer
  to \cite{enlcat}. But $\str{Et(X)}$ is in general not again
  a site, and so the objects of $\str{\etshv{X}}$ are not
  really sheaves on some site.
\end{bem}

\vspace{\abstand}

\begin{satz}
  For all $X\in Sch/S$, the canonical functor
  $*:\etshv{X}\ra\str{\etshv{\str{X}}}$ is
  exact and maps injective objects to injective
  objects.
\end{satz}

\vspace{\abstand}

\begin{proof}
  Follows from \ECref{corinj}.
\end{proof}

\vspace{\abstand}

For an $X\in\str{(Sch/S)}$, we call an object of $\str{\etshv{X}}/\str{\constr{X}}/\str{\locconst{X}}$ a
\emph{*sheaf/*con\-struc\-ti\-ble sheaf/*locally constant sheaf on $X$}.

By $Ab$ we denote the category of abelian groups in our universe.

\vspace{\abstand}

\begin{bspe}
  \mbox{}\\[-5mm]
  \begin{enumerate}
  \item Let $X$ be a scheme. For each abelian group $A\in Ab$ we get the constant sheaf $A_X$ on $X$
    (which is just the sheafification of the functor $U\mapsto A$). This gives an exact functor
    $Ab\ra\etshv{X}$, and by enlarging this, we get an exact functor
    $$\str{Ab}\ra\str{\etshv{X}},\; A\mapsto A_X.$$
  \item Let $X$ again be a scheme, and let $\G_{m,X}$ be the étale sheaf on
    $X$ which is defined by
    $$U\mapsto \Gamma(U,\cO_X)^*.$$
    If $n\in\N$ is a natural number, we have the morphism of sheaves
    $$n:\G_{m,X}\ra\G_{m,X},\; f\mapsto f^n,$$
    and if $n$ is invertible in $\Gamma(X,\cO_X)$, then this morphism
    is a surjection in $\etshv{X}$ whose kernel is denoted by $\mu_n:=(\mu_n)_X$.
    For all invertible $n\in\N$ this gives the so called \emph{Kummer sequence}
    $$0\ra (\mu_n)_X\ra \G_{m,X}\xrightarrow{n} \G_{m,X}\ra 0.$$
    If furthermore there is a primitive $n$-th root of unity $\zeta_n\in\Gamma(X,\cO)$,
    there is an isomorphism
    $$(\Z/n\Z)_X\xrightarrow{\sim}(\mu_n)_X,\; 1\mapsto \zeta_n.$$

    Now let $\str{\G_{m,X}}$ be the corresponding *sheaf on $X$. Suppose that $n\in\str\N$
    is a nonstandard natural number which is invertible in $\str{\Gamma(X,\cO_X)}$. We again get
    a morphism
    $$n:\str{\G_{m,X}}\ra\str{\G_{m,X}}$$
    which by transfer is surjective in $\str{\etshv{X}}$. We again denote by
    $\mu_n=(\mu_n)_X$ the kernel of this morphism. So for each $n\in\str{\N}$ which is invertible on $\str{X}$,
    we get a Kummer sequence
    $$0\ra(\mu_n)_X\ra\str{\G_{m,X}}\xrightarrow{n}\str{\G_{m,X}}\ra 0.$$
    Although $\str{\G_{m,X}}$ is a \emph{standard} *sheaf, the *sheaf $\mu_n$
    for $n\in(\str{\N}\setminus\N)$ is \emph{not} a standard *sheaf. Furthermore,
    if $\zeta_n$ is an $n$-th primitive root of unity in $\str{\Gamma}(X,\cO_X)$,
    we again have an isomorphism
    $$(\str{\Z}/n\str{\Z})_X\xrightarrow{\sim} (\mu_n)_X.$$
  \item For two étale sheaves $\cF$ and $\cG$ we can form the tensor product $\cF\otimes\cG$
    which is again an étale sheaf. Therefore for every two étale
    *sheaves, we also can form the tensor product with the usual properties. In fact we can do even
    more: For each *finite internal family $\{\cF_i\}_{i\in I}$ of étale *sheaves,
    we can define the tensor product $\otimes_{i\in I} \cF_i$ which is again a *sheaf.
    For example we have for each $n\in\str{\N}$ and $k\in\str{\N}$ the étale
    *sheaf $\mu_n^{\otimes k}$.
  \end{enumerate}
\end{bspe}

\vspace{\abstand}

\noindent
We also get abelian fibrations
$$\str{\text{ModConstr}(X)}\ra\str{(\finring)^{op}}$$
and
$$\str{\text{ModShv}_{et}(X)}\ra\str{(\finring)^{op}}.$$
We call the objects in $\str{\finring}$ \emph{*finite rings}. First of all
they are just objects in the category $\str{\finring}$. But
we have a canonical functor
$$\str{\finring}\ra \text{Cat. of all commutative rings}$$
which is faithful but not fully faithful. For details compare \ECref{satzrmodfun}.

\vspace{\abstand}

\begin{bsp}
  For all $n\in\str{\N}$ and all $X\in Sch/S$, the sheaves $(\mu_n)_X$ are in
  $\str{\modconstr{X}{\str{\Z}/n\str{\Z}}}.$
\end{bsp}

\vspace{\abstand}

\noindent
As before we get
\begin{satz}
  For each $S$--scheme $X\in \str{(Sch/S)}$ and each *finite ring $R\in\str{\finring}$,
  the fibre $\str{\modconstr{X}{R}}$
  is an abelian category. For $X\in Sch/S$ and $R\in\finring $ we have
  $$ \str{\modconstr{\str{X}}{\str{R}}}=\str{(\modconstr{X}{R})}$$
  and
  $$ \str{\modshv{\str{X}}{\str{R}}}=\str{(\modshv{X}{R})}.$$
\end{satz}

\vspace{\abstand}

\subsection{Étale cohomology of nonstandard schemes and sheaves}

\mbox{}\\
\vspace{\abstand}

\noindent
For all $i\in\N_0$, the $i$--th étale cohomology is the contravariant functor
$$ H^i(\_,\_) :\fib_{\etshvo} \ra  Ab, (X,\cF) \mapsto H^i_{et}(X,\cF),$$
where $H^i_{et}(X,\cF)$ is just the $i$--th right derivbed functor of the
section functor $\Gamma(X,\_):\etshv{X}\ra Ab$ applied to $\cF$.
If we enlarge this, we get for each $i\in\str{\N}$ a contravariant functor

$$\str{H}^i_{et}(\_,\_): \str{\fib_{\etshvo}}\ra\str{Ab}.$$

Here $\str{Ab}$ is first of all only an abelian category. But it is easy to
see that there is a functor from $\str{Ab}$ to the category of all abelian
groups which is faithful but not fully faithful (compare again \ECref{satzrmodfun}).
Furthermore this functor is exact by \ECrefs{satzab}{satzabex}.

\vspace{\abstand}

\noindent
For all $X\in\str{(Sch/S)}$ by proposition \ref{elestar} we can build the right derived
functor
$$\RR^i\str{H^0(X,\_)}:\str{\etshv{X}}\ra\str{Ab}.$$
By the next proposition this gives the same as the above definition.

\vspace{\abstand}

\begin{satz}
  For all $X\in\str{(Sch/S)}$, $\cF\in\str{\etshv{X}}$ and $i\in\N_0$, we have
  a canonical isomorphism
  $$ \RR^i\str{H^0(X,\_)} (\cF)\ra \str{H^i(X,\cF)}$$
\end{satz}

\vspace{\abstand}

\begin{proof}
  This follows directly from \ECref{corrflf}
\end{proof}
In particular we see that nonstandard cohomology is always a right
derived functor of an left exact functor.
\vspace{\abstand}

\noindent
For standard schemes and sheaves we have furthermore the following

\begin{satz}
  For all $i\in\N_0$, $X\in Sch/S$ and $\cF\in\str{\etshv{\str{X}}}$, the object
  $\str{H}^i_{et}(X,\cF)$ is just the $i$--th right derived functor of
  $\str{\Gamma}(X,\_):\str{\etshv{\str{X}}}\ra\str{Ab}$ applied to $\cF$.
  Furthermore  we have the following
  commutative diagram
  $$\xymatrix{  \etshv{X}  \ar[d]_{H^i_{et}(X,\_)} \ar[r]^<<<<{*}  & \str{(\etshv{X})}=\str{\etshv{\str{X}}}\ar[d]^{\str{H}^i_{et}(\str{X},\_)}   \\
    Ab  \ar[r]^{*}                                               &     \str{Ab}  }  $$
\end{satz}

\vspace{\abstand}

\begin{proof}
  This follows again from statement \ECref{corrflf}.
\end{proof}

\vspace{\abstand}

\begin{cor}
  If for an $X\in Sch/S$ and $\cF\in\etshv{X}$ the group $H^i_{et}(X,\cF)$ is finite,
  then we have an isomorphism $$H^i_{et}(X,\cF)\simeq\str{H}^i_{et}(\str{X},\str{\cF}).$$
\end{cor}

\vspace{\abstand}

\begin{proof}
  This follows from the above proposition and the fact that the enlargement of a finite set is just
  that finite set itself.
\end{proof}

\vspace{\abstand}

\subsection{Cohomology with compact support}

\mbox{}\\
\vspace{\abstand}

We are more interested in the cohomology with compact support.
For that it is most convenient to work with derived categories.
We denote by $\derconstr{X}{R}$ the full subcategory of $\cD^b(\modshv{X}{R})$
of complexes which have constructible cohomology sheaves.
Then we look at the pseudo functor
$$(Sch/S)^{compac}\times(\finring)^{op} \ra \text{Cat of triangulated categories in U}$$
that maps a pair $(X,R)$  to the triangulated category $\derconstr{X}{R}$
and a morphism $(f,\varphi):(X,R)\ra (Y,R')$ in $Sch/S\times (\finring)^{op}$
to the functor
$$(\RR f_!)_{\varphi}:\derconstr{X}{R}\ra \derconstr{Y}{R'}.$$
Here $\RR f_!$ is the direct image with proper support,
and $(\_)_\varphi$ means that we define the $R'$--module structure via the ring homomorphism
$\varphi:R'\ra R$. This defines a triangulated cofibration
$$\cD^b_{constr}\ra(Sch/S)^{compac}\times\finring^{op}.$$
By \ECref{abfib} the enlargement of this fibration
$$\str{\cD^b_{constr}}\ra\str{(Sch/S)}^{compac}\times\str{\finring^{op}}$$
is again a triangulated cofibration.
Then we go back again and look at the associated pseudo functor
$$\str{(Sch/S)^{compac}}\times\str{\finring}^{op}\ra \str{\text{(Cat. of triangulated cat. in U)}}$$
and for an object $(X,R) \in \str{(Sch/S)}\times\str{\finring}^{op}$
by $\str{\derconstr{X}{R}}$ we denote the image under this functor resp. the fibre of
the above fibration. For a morphism
 $(f,\varphi):(X,R)\ra (Y,R')$ in $\str{(Sch/S)}\times (\str{\finring})^{op}$
we denote the image of $(f,\varphi)$ under the above functor again by
$$(\RR f_!)_{\varphi}:\str{\derconstr{X}{R}}\ra \str{\derconstr{Y}{R'}}.$$

\vspace{\abstand}

In this paper we are more interested in enlarging the sheaves and not so much
in enlarging the schemes.
Therefore we now only consider --- for each $X\in Sch/S$ --- the cofibrations of
triangulated categories
$$\derconstrs{X}\ra (\finring)^{op}$$
and
$$\cD^b(\modshvs{X})\ra (\finring)^{op}$$
and the enlargements
$$\str{\derconstrs{X}}\ra \str{(\finring)}^{op}$$
and
$$\str{\cD^b}(\modshvs{X})\ra \str{(\finring)^{op}}$$

\vspace{\abstand}

The notion of constructibility in $\modshv{X}{R}$ for all $R\in\finring$ gives
us a notion of *constructibility in $\str{\modshv{X}{R}}$ for
all $R\in\str{\finring}$. From the transfer principle we get the following

\vspace{\abstand}

\begin{satz}
  Let $(\str{\cD^b(\modshvs{X})}_{constr}$ be the full subcategory of
  $\str{\cD^b(\modshvs{X})}$ of those *complexes which have
  *constructible cohomology sheaves (for all $i\in\str{\Z}$ !).
  Then there is a canonical isomorphism
  $$(\str{\cD^b(\modshvs{X}))}_{constr} \xrightarrow{\sim} \str{\cD^b_{constr}(X)}.$$
\end{satz}

\vspace{\abstand}

\begin{proof}

\end{proof}

\vspace{\abstand}

Let $\str{\derconstrsf{X}}$ resp. $\cD^{fb}(\modshvs{X})$ be the full subcategory
of those *complexes which are cohomologically bounded by a standard (!)
natural number. For more detail we refer to \cite[sect. 7]{enlcat}. With that we get

\vspace{\abstand}

\begin{satz}\label{derunder}
  For all $X\in Sch/S$ and each $R\in\str{\finring}$,
  we have a canonical isomorphism of triangulated categories
  $$\str{\cD^{fb}}(\modshvs{X})(R)\xrightarrow{\sim}\cD^b(\str{\modshvs{X}}(R))$$
  and
  $$\str{\derconstrf{X}{R}}\ra\cD^b(\modconstr{X}{R}).$$
\end{satz}

\vspace{\abstand}

\begin{proof}
  This follows from \ECref{fibversion}
  and from the fact that for all $R\in\finring$
  $\derconstr{X}{R}=\cD^b(\modconstr{X}{R})$ (compare \cite[prop.4.6, page 93]{SGA4.5}).
\end{proof}

\vspace{\abstand}

\begin{bem}
  In fact that even gives an isomorphism of cofibred categories.
\end{bem}

\vspace{\abstand}

\begin{bem}
We are mainly interested in $\derconstrf{X}{\str{\Z}/P\str{\Z}}$ where
$P$ is an infinite prime in $\str{\Z}$.
From an external point of view,
the ring $\str{\Z}/P\str{\Z}$ is a field of characteristic $0$,
which leads us to think that $\derconstrf{X}{\str{\Z}/P\str{\Z}}$ is a good
replacement for the derived category of $\Q_l$--sheaves on $X$ (for a standard prime $l$).
One advantage of $\derconstrf{X}{\str{\Z}/P\str{\Z}}$ is that it
really is the derived category of an abelian category.
\end{bem}

\vspace{\abstand}

\noindent
We further get the following property for the direct image with proper support:

\vspace{\abstand}

\begin{satz}\label{sternderivation}
  Let $X\ra Y$ be a compactifiable morphism in $Sch/S$,
  and let $f=\bar f\circ j$ be a factorization of $f$
  with an open immersion $j:X\hookrightarrow\bar X$
  and a proper morphism $\bar f:\bar X \ra Y$.
  Then for each $R\in\finring$ we get the following commutative diagram:
  $$
  \xymatrix@C=25mm{
    \cD^b(\str{\modconstr{X}{R}}) \ar[r]^{\RR \str{(\bar f_*)}\circ \str{(j_!)}} \ar[d]^{\wr}   &
    \cD^b(\str{\modconstr{Y}{R}}) \ar[d]_{\wr} \\
    \str{\derconstrf{X}{R}} \ar[r]^{\str{(\RR f_!)}}   &  \str{\derconstrf{Y}{R}}.
    }
  $$
\end{satz}

\vspace{\abstand}

\begin{proof}
  For all $R\in\finring$, the cohomological dimension of the functor
  $$f_!:\modconstr{X}{R}\ra\modconstr{Y}{R}$$
  is less or equal to the relative dimension of the morphism $f$ (which is a \emph{finite} number).
  Therefore for all $R\in\str{\finring}$, the same is true for the functor
  $$\str{f_!}:\str{\modconstr{X}{R}}\ra\str{\modconstr{Y}{R}}.$$
  This shows that $\str{\RR f_!}$ restricts to the finitely bounded complexes,
  and then the proposition follows from \ECref{fibversion}.
\end{proof}

\vspace{\abstand}

This proposition justifies the following notation: In the situation
of the proposition, for all $R\in\str{\finring}$,
we denote the functor $$\str{(\RR f_!)}: \str{\derconstrf{X}{R}}\ra\str{\derconstrf{Y}{R}}$$
again by $\RR f_!$,
and for $i\in\Z$ we set
$$\RR^if_!:=h^i\circ \RR f_!,$$
where $h^i$ is the $i$-th cohomology functor
$$h^i:\str{\derconstrf{X}{R}}\ra\str{\modconstr{Y}{R}}.$$

\vspace{\abstand}

\noindent
For two composable compactifiable morphisms we have the following result:

\vspace{\abstand}

\begin{satz}
  Let $X,Y,Z\in (Sch/S)$,
  let $f:X\ra Y$ and $g:Y\ra Z$ be two compactifiable morphisms in $Sch/S$,
  and let $R\in\finring$.
  Then we have the following commutative diagram of triangulated categories:
  $$
  \xymatrix{
    \str{\derconstrf{X}{R}} \ar[rd]^{\RR f_!} \ar[dd]^{\RR (g\circ f)_!}   &    \\
                                                          & \str{\derconstrf{Y}{R}} \ar[dl]^{\RR g_!} \\
    \str{\derconstrf{Z}{R}}
  }
  $$
  In particular, for $\cF\in\str{\modconstr{X}{R}}$ we have a spectral sequence
  $$ E^{p,q}_2=\RR^p g_!(\RR^q f_!(\cF)) \Rightarrow \RR^{p+q} (g\circ f)_!(\cF).$$
\end{satz}

\vspace{\abstand}

\begin{proof}
  This follows from the equality
  $\RR g_!\circ \RR f_!=\RR(g\circ f)_!$
  in the standard world and from proposition \ref{sternderivation}.
\end{proof}

\vspace{\abstand}

\noindent
For a ring $R$, let $\Mod{R}$ be the category of $R$--modules in our universe.
The pseudo functor
$$\finring\ra \text{Cat. of triang. categories in U}$$
which sends $R$ to $\cD^b(\Mod{R})$ defines --- as above --- a cofibration
$$\cD^b(\text{Mod})\ra\finring^{op},$$
and we again get the enlargement
$$\str{\cD^b(\text{Mod})}\ra\str{\finring^{op}}$$
and the analogue of proposition \ref{derunder}.

\vspace{\abstand}

\begin{satz}
  For each $R\in\str{\finring}$ there is a canonical isomorphism of triangulated
  categories
  $$\str{\cD^{fb}}(\text{Mod})(R)\xrightarrow{\sim}\cD^b(\str{\Mod{R}}).$$
\end{satz}

\vspace{\abstand}

\begin{proof}
  This can be proven in the same way as proposition \ref{derunder}.
\end{proof}

\vspace{\abstand}

Now we assume that $S$ has finite cohomological dimension
in the étale topology. Then the derived functor of the global section
functor
$$\Gamma(S,\_):\etshv{X}\ra Ab$$
defines for each $R\in\finring$ a functor
$$\derconstr{S}{R}\ra\cD^b(\Mod{R})$$ and so also for each $R\in\str{\finring}$ a functor
$$\str{\derconstrf{S}{R}}\ra\str{\cD^{fb}}(\text{Mod})(R).$$
In the same way as above we get the analogue of proposition \ref{sternderivation}.

\vspace{\abstand}

\begin{satz}
  Let $S$ be a scheme of finite cohomological dimension in the étale topology.
  Then we get for each $R\in\str{\finring}$ the following commutative diagram
  $$
  \xymatrix@C=25mm{
    \cD^b(\str{\modconstr{S}{R}})\ar[d]^{\wr} \ar[r]^-{\RR \str{\Gamma}(S,\_)} & \cD^b(\str{\Mod{R}})\ar[d]_{\wr}  \\
    \str{\derconstrf{S}{R}} \ar[r]^{\str{(\RR\Gamma(S,\_))}}    & {\str{\cD^{fb}}(\text{Mod})(R).}
  }
  $$
\end{satz}

\vspace{\abstand}

\begin{proof}
  This can be proven in the same way as proposition \ref{sternderivation}.
\end{proof}

\vspace{\abstand}

\noindent
We will again omit the * in the notation.\\

Now if $f:X\ra S$ is a compactifiable morphism,
if $i\in\str{\N_0}$,
if $R\in\str{\finring}$,
and if $\cF\in\str{\modconstr{X}{R}}$,
then we call the object
$$H^i_c(X,\cF):=\RR^i\Gamma(S,\RR f_!\cF)$$
in $\str{\Mod{R}}$ the \emph{$i$-th cohomology of $\cF$ on $X$ with proper (or compact) support}.
As already mentioned several times before,
$H^i_c(X,\cF)$ can also be seen as a module over the ring $R$.

Now we obviously get the following spectral sequence:

\vspace{\abstand}

\begin{satz}\label{hochserre}
  Let $X\ra S$ be a compactifiable morphism,
  let $i\in\str{\N_0}$,
  let $R\in\str{\finring}$,
  and let $\cF\in\str{\modconstr{X}{R}}$.
  Then we have a spectral sequence
  $$E_2^{p,q}:=H_c^i(S,\RR^p f_!\cF)\Rightarrow H_c^{p+q}(X,\cF).$$
\end{satz}

\vspace{\abstand}

\begin{proof}
  This is clear by the previous propositions.
\end{proof}

\vspace{\abstand}


\section{Nonstandard cohomology of schemes over a field}

\mbox{}\\
\vspace{\abstand}

In this chapter we want to apply the results of the first chapter
to the case of varieties, i.e. schemes which are separated and of finite type over a field.
The category of *constructible sheaves
of $K$--modules for a field $K\in\str{\finring}$ on the spectrum of a separably closed field is equivalent
to the category of *finite dimensional internal vector spaces
over the field $K$, and such spaces are in general not finite dimensional.
However, in this chapter we want to show
that at least the geometric cohomology (see below)
of varieties
with coefficients in finite dimensional constant sheaves
is again finite dimensional.

\vspace{\abstand}

\subsection{Étale sheaves over a field}

\mbox{}\\
\vspace{\abstand}

\def\contrep{\text{ContRepr}(Gal(k_s/k))}
\def\gal{Gal(k_s/k)}

Let $k$ be a field, let $k_s$ be a separable closure of $k$, and
let $R$ be a finite ring.
We denote by $\contrep(R)$ the category of continuous $\gal$--$R$--modules
in our universe $U$.
Then the functor
$$F\mapsto \varinjlim_{\text{\tiny $k'/k$ finite separable}}F(k')$$
gives us an equivalence of categories between $\modshv{\spec{k}}{R}$
and $\contrep(R)$.
The enlargement of this gives us for each $R\in\str{\finring}$
an equivalence between $\str{\modshv{\spec{k}}{R}}$ and
$\str{\contrep(R)}$. By \ECref{satzrmodfun} the category $\str{\contrep(R)}$
can be interpreted as specific (i.e. *continuous and internal)
$R$--representations of the group $\str{\gal}$.
If $k$ is a finite field of characteristic $p$, and if $F_p\in\gal$ is the Frobenius,
then $F_p$ also acts on each object of $\str{\contrep(R)}$.

\vspace{\abstand}

\noindent
We denote by
$$H^i(\str{\gal},\_): \str{\modshv{\spec{k}}{R}} \ra \str{\Mod{R} }$$
the functor $\RR^i\Gamma(\spec{k},\_)$, and we call this functor
\emph{nonstandard Galois cohomology}.

\vspace{\abstand}
Now let $X$ be a variety over a field $k$ with structure morphism $f:X\ra\spec{k}$.
Let $i\in\str{\N_0}$,
let $R\in\str{\finring}$ be a *finite ring, and let $\cF$  be a *constructible sheaf of $R$--modules.

Then we call the objects
$$H_c^i(\bar{X},\cF):=\RR^i f_! \cF\in\str{\contrep(R)}$$
the \emph{geometric cohomology of $X$ with coefficients in $\cF$},
and we call
$$H^i_c(X,\cF):=H^i(\str{\gal},\RR f_!\cF)\in\str{\Mod{R}}$$
the \emph{arithmetic  cohomology of $X$ with coefficients in $\cF$}.
One can easily see that $H^i_c(\bar X,\_)$ is the enlargement
of the functor $\cF\mapsto H^i_c(X\otimes_k k_s,\cF_{k_s})$.

\vspace{\abstand}

\noindent
From proposition \ref{hochserre} we immediately get the following Hochschild-Serre
spectral sequence:
\vspace{\abstand}

\begin{satz}
  In the above situation we get a spectral sequence
  $$E_2^{p,q}:=H^p(\str{\gal},H^q_c(\bar X,\cF))\Rightarrow H^{p+q}_c(X,\cF)$$
\end{satz}
\begin{proof}
  This is a special case of \ref{hochserre}
\end{proof}

\vspace{\abstand}

\noindent
By the transfer principle we see that for a *constructible sheaf $\cF$ on $X$ we have
$$H^i_c(\bar X,\cF)=0 \text{ for } i>2\cdot dim\,X$$
and
$$H^i_c(X,\cF)=0 \text{ for } i>2\cdot dim\,X + cohdim\,k.$$

\vspace{\abstand}

Now let $P\in\str{\N}$ be an \emph{infinite} prime,
and set $R:=\str{\Z}/P$.
Until now we only know that $H^i_c(\bar X,\str{\Z}/P)$ is *finite dimensional.
We want to show that it actually is \emph{finite} dimensional.

\begin{satz}\label{finite}
  In the above situation, the $\str{\Z}/P$--vector spaces $H^i_c(\bar X,\str{\Z}/P)$ are finite dimensional.
\end{satz}

\vspace{\abstand}

\begin{proof}
  This follows from the transfer principle and the next proposition.
\end{proof}

\vspace{\abstand}

\begin{satz}
  Let $X$ be a variety over an algebraically closed field.
  Then there is a constant $C\in\N$ with the property that for all
  (standard) primes $l\in\N$ and all $i\in\N_0$ the $\Z/l$--rank of $H^i_c(X,\Z/l)$
  is less or equal $C$.
\end{satz}

\vspace{\abstand}

\begin{proof}
  By de Jong's alteration we can assume that $X$ is smooth. Further it follows easily that we can assume
  that $X$ is projective. From \cite{gabber} we know that the torsion
  in $H^i(X,\Z_l)$ (l-adic cohomology) is zero for almost all $l$. By the independence of
  $l$ for smooth and projective Varieties it follows that the rank of $H^i(X,\Z_l)$ is the
  same for almost all $l$. Finally, we have the long exact sequence
  $$\dots\ra H^i(X,\Z_l)\xrightarrow{\cdot l} H^i(X,\Z_l)\ra H^i(X,\Z/l) \ra H^{i+1}(X,\Z_l)\xrightarrow{\cdot l}
                                 H^{i+1}(X,\Z_l)\ra\dots$$
  From this the proposition follows.
\end{proof}

\vspace{\abstand}

\begin{bem}
  The proof of proposition \ref{finite} shows that the finiteness of
  the dimension of $\str{H}^i_c(\bar{X},\str\Z/P)$ for
  $P\in\str{P}\setminus P$ seems to be a much deeper result than the
  finiteness of ordinary $l$-adic cohomology. Somehow the cohomology
  groups $\str{H}^i_c(\bar{X},\str\Z/P)$ contains the $l$-adic
  cohomology for infinitely many standard primes $l$.

\end{bem}

\vspace{\abstand}

\subsection{Poincaré duality}

\mbox{}\\
\vspace{\abstand}

First we want to recall what the Poincaré duality in the standard world means.
For that let $X$ be a smooth connected scheme of dimension $d$ over a separably closed
field. Then there is for each $n\in\Z$ a canonical trace map
$$H_c^{2d}(X,\mu_n^{\otimes d})\ra \Z/n$$
which is an isomorphism.

Furthermore for each constructible $\Z/n$--module
$\cG$ the cup product pairing
$$Ext^p(\cG,\mu_{n,X}^{\otimes d})\times H_c^{2d-p}(X,\cG)\ra H_c^{2d}(X,\mu_n^{\otimes d})\simeq \Z/n$$
is non degenerate.

Now we want to translate this statement to our enlargement.
So let $n\in\str{\N}$ be a hyper natural number.
First of all we remark that for two $\str{\Z}/n$--module *sheaves $\cF$ and $\cG$ on $X$,
we can identify the object $\str{Ext}^p(\cF,\cG)$
with the $p$--th $Ext$ group in the category $\str{\modshv{X}{\str{\Z}/n}}$ for all $p\in\N_0$,
and it is easy to see that the transfer of the cup product pairing is the cup product pairing
in $\modshv{X}{\str{\Z/n}}$.
Then by the transfer principle we get

\vspace{\abstand}

\begin{satz}
  Let $X$ be a smooth connected scheme of dimension $d$ over a separably closed field,
  and let $n\in\str{\Z}$.
  Then we get a canonical trace map
  $$H_c^{2d}(X,\mu_n^{\otimes d})\ra \str{\Z}/n,$$
  which is an isomorphism. Furthermore for each *sheaf $\cF\in\modconstr{X}{\str{\Z/n}}$
  the cup product pairing
  $$Ext^p(\cG,\mu_{n,X}^{\otimes d})\times H_c^{2d-p}(X,\cG)\ra H_c^{2d}(X,\mu_n^{\otimes d})\simeq \str{\Z}/n$$
  is non degenerate.
\end{satz}

\vspace{\abstand}

\begin{proof}
  The only thing we have to add to our above remarks is the fact that the enlargement
  of a non degenerate pairing is again a non degenerate pairing. But this
  follows easily from the transfer principle.
\end{proof}

\vspace{\abstand}

\subsection{The cycle map}

\mbox{}\\
\vspace{\abstand}

Let $X$ be a smooth variety over an algebraically closed field. We denote by
$Z^r(X)$ the group of algebraic cycles of codimension $r$ on $X$. Further we
denote by $CH^r(X)$ the Chow--group of codimension--$r$--cycles.

\vspace{\abstand}

For each prime $l\in\N$ which is invertible on $X$ and each $r\in\N_0$ there is
a cycle map
$$cl_X: Z^r(X)\rightarrow H^{2r}(X,\mu_l^{\otimes r}).$$

We write $H^*(X,\Z/l)$ for $\oplus_r H^{2r}(X,\mu_l^{\otimes r})$.
With the cup product as multiplication, $H^*(X,\Z/l)$ becomes an anti commutative graded ring.

\vspace{\abstand}

Now, if $(dim(X)-1)!$ is invertible in $\Z/l$,
the cycle map induces a morphism of graded rings
$$cl_X:CH^*(X)\rightarrow H^*(X,\Z/l)$$
which is functorial in $X$ (compare for example \cite[VI.\S10]{milne}).

\vspace{\abstand}

Now let $P\in\str{\N}$ be an infinite prime which is invertible on $X$,
and let $r\in\N_0$.
By transfer we again get a map
$$cl_X:\str{(Z^r(X))}\rightarrow H^{2r}(X,\mu_P^{\otimes r}),$$
and as for finite primes,
the cup product again induces the structure of an anti commutative graded ring on
$$H^*(X,\str{\Z}/P):=\oplus_r H^{2r}(X,\mu_P^{\otimes r}).$$
Because $dim(X)$ is finite and $P$ is infinite,
$(dim(X)-1)!$ is always invertible in $\str{\Z}/P$,
and we always get a morphism of graded rings
$$\str{(CH^r(X))}\rightarrow H^*(X,\str{\Z}/P).$$
If we compose this with the canonical morphism
$CH^*(X)\rightarrow \str{(CH^*(X))}$,
we again get a morphism of anti commutative graded rings
$$cl_X: CH^*(X) \rightarrow H^*(X,\str{\Z}/P)$$
which is functorial in $X$.

\vspace{\abstand}

\subsection{Further properties}

\mbox{}\\
\vspace{\abstand}

We would like to conclude this section by a general and quite unprecise
statement about the other fundamental properties of étale cohomology. We state
this in the following

\begin{bem}[Meta-theorem]
   Each appropriate theorem about \'etale cohomology with finite coefficients has
   its nonstandard version for *finite coefficients.
\end{bem}

\vspace{\abstand}

Above, we gave the examples of Poincaré duality and of the existence of cycle maps.
Other examples are purity, Künneth formula, base change theorems, \dots .

\vspace{\abstand}

\begin{thm}
Let $P$ be an infinite prime
and let $k$ be an algebraically closed field of finite characteristic prime
to $P$.
The contravariant functor
$$X\mapsto\bigoplus_{n\geq 0}H^n(X,\str{\Z}/P)$$
from the category of smooth projective varieties over k
to the category of anti commutative graded $\str{\Z/P}$--algebras,
is a Weil cohomology theory in the sense of \cite{kl}.
\end{thm}

\vspace{\abstand}

\begin{proof}
The fact that the cohomology groups are finite dimensional is just
proposition \ref{finite}. Except the hard Lefschetz theorem everything
else follows from the transfer principle and the corresponding
statement for $\Z/l$ cohomology. The hard Lefschetz follows from
\cite{katzmessing}.
\end{proof}

\begin{bem}
Toma\v{s}i\'{c} constructed in \cite{tomasic} in a similar way a Weil cohomology
theory for smooth projective varieties. His approach is more direct
and less technical
where our approach allows more flexibility with the coefficients and
also exhibits this cohomology theory as a derived functor.
\end{bem}

\vspace{\abstand}


\section{Comparison with singular cohomology}

\mbox{}\\
\vspace{\abstand}

In this sections we consider smooth and projective schemes over $\C$. We want to compare the
nonstandard étale cohomology with constant *finite coefficients with singular
cohomology of the associated complex analytic space. For the case of finite
coefficients we refer to standard books on étale cohomology.

\vspace{\abstand}

We denote by $\finab$ the category of finite abelian groups.
For each topological space $X$ in our universe we have the
functors
$$H^i(X,\_):\finab\ra Ab$$
of singular cohomology. This gives us functors
$$\str{H}^i(X,\_):\str{\finab}\ra \str{Ab}.$$
Each object of $\str{\finab}$ and $\str{Ab}$ defines
an abelian group (compare \ECref{satzab}), and we will use the same name for
this abelian group.
The following proposition tells us that we do not get
anything new in this case.

\vspace{\abstand}

\begin{satz}\label{singular}
  For a topological space $X$,
  an $i\in\N_0$
  and a *finite abelian group $M$  there is a canonical isomorphism
  $$\str{H}^i(X,M)\xrightarrow{\sim} H^i(X,M)$$.
\end{satz}

\begin{proof}
  Consider the cohomology theory which maps a topological space $X$
  to $\str{H}^i(X,M)$. From the transfer principle it follows that
  this cohomology theory satisfies the axioms of Eilenberg and Steenrod.
  Further the  dimension axiom
  $$\str{H}^i(\{pt\},M)=0 \text{ for } i\not=0$$
  holds. Because $\str{H}^0(X,M)=M$, the theorem follows from
  the universality of singular cohomology (compare for example
  \cite[Th. 4.59]{hatcherat}).
\end{proof}

\vspace{\abstand}

Now let $X$ be a smooth projective scheme  over $\C$. We denote by
$X(\C)$ the associated complex analytic space. Then we have
the following comparison theorem:

\begin{satz}
  Let $X$ be a smooth and projective scheme over $\C$,
  and let $M\in\str{\bigl(\finab\bigr)}$.
  Then there is a canonical isomorphism
  $$H^i(X(\C),M)\xrightarrow{\sim} H^i(\str{X},M_{\str{X}}).$$
\end{satz}

\vspace{\abstand}

\begin{proof}
  It is known that for all $M\in\finab$ there is a canonical
  isomorphism
  $$H^i(X(\C),M)\ra H^i(X,M_X).$$
  By transfer it follows that for each $M\in\str{\bigl(\finab\bigr)}$ there
  is a canonical isomorphism
  $$\str{H}^i(X(\C),M)\ra H^i(X,M_X).$$
  Now the proposition follows with proposition \ref{singular}.
\end{proof}

\vspace{\abstand}


\def\F{{\mathbb F}}
\def\R{{\mathbb R}}
\def\C{{\mathbb C}}
\def\N{{\mathbb N}}
\def\Nn{{{\mathbb N}_0}}
\def\Np{{{\mathbb N}_+}}
\def\Q{{\mathbb Q}}
\def\Z{{\mathbb Z}}
\def\P{{\mathbb P}}
\def\Ns{\str{\mathbb N}}
\def\Nns{\str{{\mathbb N}_0}}
\def\Nps{\str{{\mathbb N}_+}}
\def\Cs{\str{\mathbb C}}
\def\Ps{\str{\mathbb P}}
\def\Qs{\str{\mathbb Q}}
\def\Rs{\str{\mathbb R}}
\def\Zs{\str{\mathbb Z}}
\def\Qab{{\Q^{\text{ab}}}}
\def\Qb{\bar{\Q}}
\def\ordp{{\text{ord}_p}}
\def\p{{\mathfrak p}}
\def\pinf{{(p^\infty)}}
\def\resp{resp.\ }
\def\Ob{\mathrm{Ob}}
\def\cA{\mathcal A}
\def\cB{\mathcal B}
\def\cC{\mathcal C}
\def\cD{\mathcal D}
\def\cI{\mathcal I}
\def\cJ{\mathcal J}
\def\cL{\mathcal L}
\def\cP{\mathcal P}
\def\cR{\mathcal R}

\def\Ensfin{{\mathrm{Ens}}^{\mathrm{fin}}}
\def\finG{\fin{\tilde{G}}}
\def\strfin{$\str{\text{finite}}$}
\newcommand\MorC[1]{{\mathrm{Mor}}_{#1}}
\newcommand\Modfin[1]{{\mathrm{Mod}_{#1}^{\mathrm{fin}}}}
\def\image{\text{im}\,}
\def\kernel{\text{ker}\,}
\def\cokernel{\text{coker}\,}
\newcommand\rmod[1]{(#1-\underline{\mathrm{Mod}})}
\def\mR{\mathfrak{R}}
\def\Ensigns{\underline{\mathrm{Ens}}}

\newcommand\pr[1]{\mathrm{Proj}\,(#1)}
\newcommand\artinrees[1]{\mathrm{AR}\,(#1)}
\newcommand\artinreesl[1]{\mathrm{AR}_l(#1)}
\def\Abg{\mathrm{Ab}}
\newcommand\arfun[2]{\Upsilon_{#1}#2}
\newcommand\arstrfun[1]{\Xi #1}
\newcommand\arzlhfun[2]{\Psi_{#1}#2}

\section{Comparison with $l$-adic cohomology}

\vspace{\abstand}

Let $l$ be a (standard) prime number, let $h\in\str{\N}$ be an infinite natural number,
and let $X$ be a noetherian scheme.

In this chapter, we want to show that the category of *constructible $\str{\Z}/l^h$--modules
can serve as a substitute for the category of $l$-adic sheaves or Artin-Rees-$l$-adic sheaves on $X$
in a sense to be made precise later.

In particular, we will show that geometric $l$-adic cohomology of a variety $X$
can be computed using geometric cohomology of $X$
with coefficients in the *constructible sheaves $\str{\Z}/l^{d_1+d_2}$ and $\str{\Z}/l^{d_1}$ for two
infinite natural numbers $d_1$ and $d_2$.
For smooth and projective varieties $X$
and for almost all primes $l$,
there is an even simpler isomorphism $H^i_c(\bar{X},\str{\Z}/l^h)\otimes_{\str{\Z}}\Z_l\cong H^i_c(\bar{X},\Z_l)$.

\vspace{\abstand}

The basic idea is as follows:
If $(\cF_n)$ is an $l$-adic sheaf on $X$,
we can consider the associated *constructible $\str{\Z}/l^h$-sheaf $\str{\cF}_h$
and from this (essentially) recover $\cF$ because of $\str{\cF_n}=\str{\cF}_h/l^{n+1}$ for all $n$
(compare example \ECrefs{bsplim}{bspprojlim}).
In this way we get a faithful functor $(\cF_n)\mapsto\str{\cF}_h$
from $l$-adic sheaves to *constructible $\str{\Z}/l^h$-sheaves which reflects isomorphisms.

Unfortunately, the higher direct images of $l$-adic sheaves need not be $l$-adic again, but merely
Artin-Rees $l$-adic, i.e. isomorphic to an $l$-adic sheaf in the
Artin-Rees category\footnote{See \ref{defiar} for a definition of this term!}.
It is no longer possible to recover an Artin-Rees $l$-adic sheaf $\cF$ from $\str{\cF}_h$ as above,
but it is well known how to construct an $l$-adic sheaf which is isomorphic to $\cF$:
If $\cF'$ denotes the subsystem of $\cF$ formed by the stable images, then there is a natural number $r\in\N_0$
such that $\cG:=(\cF'_{n+r}/l^{n+1})_n$ is an $l$-adic sheaf isomorphic to $\cF$ (compare \cite[V.3.2.3]{sga5}).
But $r$ depends on $\cF$, so that this construction can not be used to define a "nice"
functor from Artin-Rees $l$-adic sheaves
to $l$-adic sheaves.

Luckily, as we have infinite natural numbers at our disposal, we are in a much more fortunate position:
First of all, we can replace the stable image by an "infinite" image
$(\image[{\str{\cF}_{n+d_1}\rightarrow\str{\cF_n}]})_n$
for an infinite number $d_2$;
second, instead of an $r$ which depends on $\cF$, we can use another infinite number $d_1$,
and finally, using our original idea for $l$-adic systems, we only have to compute $\str{\cG}_h$ to recover $\str{\cG}_n$
for finite numbers $n$.
This means we only have to consider the *constructible $\str{\Z}/l^h$-sheaf
$\image[\str{\cF}_{h+d_1+d_2}\rightarrow\str{\cF}_{h+d_1}]/l^h$ without losing essential information,
and it means we can turn this construction into a faithful functor from Artin-Rees $l$-adic systems to
*constructible $\str{\Z}/l^h$-sheaves which reflects isomorphisms,
thus enabling us to replace Artin-Rees $l$-adic sheaves by the conceptually simpler *constructible sheaves
in a functorial way.

\vspace{\abstand}

\subsection{Generalities on Artin-Rees-$l$-adic systems}

\mbox{}\\
\vspace{\abstand}

\noindent
For details and proofs, we refer the reader to \cite{sga5}[V]. ---
Let $\cA$ be an abelian category, and let $l$ be a prime number.

\vspace{\abstand}

\begin{defi}\label{defiar}
  \mbox{}\\[-5 mm]
  \begin{enumerate}
    \item
      By $\pr{\cA}$ we denote the abelian category whose objects are projective systems
      $(F_{n+1}\xrightarrow{u_{n+1}}F_n)_{n\in\N_0}$ in $\cA$
      and whose morphisms are morphisms of projective systems.
    \item
      If $F=(F_n,u_n)$ is an object of $\pr{\cA}$ and $r\in\N_0$,
      we define the \emph{shifted} system $F[r]$ as the projective system
      $(F_{n+r},u_{n+r})_{n\in\N_0}$,
      and if $f=(f_n):F\rightarrow G$ is a morphism in $\pr{\cA}$, we define the \emph{shifted} morphism
      $f[r]:F[r]\rightarrow G[r]$ by $f[r]:=(f_{n+r})_{n\in\N_0}$.
      Shifting by $r$ obviously defines an endofunctor of $\pr{\cA}$.
      Note that the $(u_n)$ induce a natural morphism
      $F[r]\rightarrow F$ in $\pr{\cA}$.
    \item
      An object $F$ of $\pr{\cA}$ is called a \emph{zero system} if there is an $r\geq 0$ such that
      the natural morphism $F[r]\rightarrow F$ is zero.
    \item\label{defiladic}
      A system $F=(F_n,u_n)$ of $\pr{\cA}$ is called \emph{$l$-adic}
      if $l^{n+1}F_n=0$ for all $n$
      and if $u_{n+1}$ induces an isomorphism $F_{n+1}/l^{n+1}F_{n+1}\xrightarrow{\sim}F_n$ for all $n$.
    \item
      By $\artinrees{\cA}$ we denote the category whose objects are the same as those of $\pr{\cA}$
      but whose morphisms are defined as
      \[
        \Mor{\artinrees{\cA}}{F}{G}:=\varinjlim_{r\geq 0}\ \Mor{\pr{\cA}}{F[r]}{G}
      \]
      with the obvious addition and composition laws.
      (Equivalently, $\artinrees{\cA}$ can be defined as the quotient of $\pr{\cA}$
      by the Serre subcategory of zero systems.)
    \item
      An \emph{Artin-Rees $l$-adic system} is an object of $\pr{\cA}$ which is isomorphic
      to an $l$-adic system \emph{in $\artinrees{\cA}$}. We denote the full subcategory of $\artinrees{\cA}$
      consisting of Artin-Rees $l$-adic systems by $\artinreesl{\cA}$.
  \end{enumerate}
\end{defi}

\vspace{\abstand}

\noindent
The straightforward proof of the following proposition and its two corollaries can be found in
\cite[V]{sga5}.

\begin{satz}\label{satzprojar}
  Let $F$, $L_1$, $L_2$ and $N$ be objects of $\pr{\cA}$ with $L_1$, $L_2$
  $l$-adic and $N$ a zero system.
  \begin{enumerate}
    \item\label{satzprojabelsch}
      $\pr{\cA}$ and $\artinrees{\cA}$ are abelian categories.
    \item\label{lemmamorphs}
      $\Mor{\pr{\cA}}{L_1}{L_2}=\Mor{\artinrees{\cA}}{L_1}{L_2}.$
    \item\label{lemmaarnull}
      $F$ is isomorphic to zero in $\artinrees{\cA}$ iff
      $F$ is a zero system.
    \item\label{lemmaladicepi}
      Let $L_1\twoheadrightarrow N$ be an epimorphism in $\pr{\cA}$. Then $N=0$.
  \end{enumerate}
\end{satz}

\vspace{\abstand}

\begin{bem}
  Note that in general, the category $\artinreesl{\cA}$ will not be abelian.
  But it is a well known fact that $\artinreesl{\cA}$ is abelian in the case where
  $\cA$ is the category of finite abelian groups
  or, more generally, the category of constructible étale sheaves on a scheme on which $l$ is invertible
  (\cite[12.11]{kiehl1}).
\end{bem}

\vspace{\abstand}

\begin{cor}\label{corarladic}
  A projective system $F\in\Ob(\pr{\cA})$ is Artin-Rees-$l$-adic if and only if
  there is an $l$-adic system $G$,
  an integer $r\geq 0$
  and an epimorphism $F[r]\twoheadrightarrow G$ \emph{in $\pr{\cA}$} whose kernel is a zero system.
\end{cor}

\vspace{\abstand}

\begin{cor}\label{corarml}
  Let $F$ be an Artin-Rees-$l$-adic system.
  Then there is an integer $s\geq 0$ such that
  $\image\bigl(F[s+n]\rightarrow F\bigr)=\image\bigl(F[s]\rightarrow F\bigr)$
  for all $n\geq 0$, i.e. $F$ has the \emph{Artin-Rees-Mittag-Leffler property}.
\end{cor}

\vspace{\abstand}

\begin{lemma}\label{lemmaarker}
  Let $0\rightarrow N\rightarrow F\rightarrow G\rightarrow 0$ be an exact sequence in $\pr{\cA}$
  with a zero system $N$ and an $l$-adic system $G$.
  Then if $N[r]\rightarrow N$ is zero, we have
  \[
    \forall m,n\geq 0:
    \image\Bigl[\kernel\bigl(F_{r+m+n}\rightarrow F_n\bigr)\rightarrow F_{m+n}\Bigr]
    \hookrightarrow l^{n+1}\cdot F_{m+n}.
  \]
\end{lemma}

\vspace{\abstand}

\begin{proof}
  We have the following commutative diagram with exact rows in $\cA$:
  \[
    \xymatrix{
      0 \ar[r] & {N_{r+m+n}} \ar[r] \ar[d]_{0} & {F_{r+m+n}} \ar[r]^{f} \ar[d]_{u} &
        {G_{r+m+n}} \ar[r] \ar@{->>}[d] & 0 \\
      0 \ar[r] & {N_{m+n}} \ar[r] \ar[d] & {F_{m+n}} \ar[r] \ar[d]_{v} & {G_{m+n}} \ar[r] \ar@{->>}[d] & 0 \\
      0 \ar[r] & {N_n} \ar[r] & {F_n} \ar[r] & {G_n\cong G_{r+m+n}/l^{n+1}} \ar[r] & 0. \\
    }
  \]
  For reasons of simplicity, using Mitchell's imbedding theorem \cite{mitchell},
  we can assume that $\cA$ is a full exact abelian subcategory of
  the category of modules over a ring. If
  $x\in\kernel(vu)$, then $f(x)\in l^{n+1}\cdot G_{r+m+n}$, i.e. $f(x)=l^{n+1}\cdot f(y)$ for a suitable
  $y\in F_{r+m+n}$ (because $f$ is surjective).
  This means $x-l^{n+1}\cdot y$ comes from $N_{r+m+n}$ and accordingly is mapped to zero in $F_{m+n}$,
  i.e. $u(x)=l^{n+1}\cdot u(y)\in l^{n+1}\cdot F_{m+n}$ as claimed.
\end{proof}

\vspace{\abstand}

\subsection{Enlargements and Artin-Rees-$l$-adic systems}

\mbox{}\\
\vspace{\abstand}

Let $*:\hat{S}\rightarrow\widehat{\str{S}}$ be an enlargement,
let $\cA$ be an $\hat{S}$-small abelian category,
let $l\in\N$ be a prime number,
and let $h\in\str{\N}\setminus\N$ be an \emph{infinite} natural number.

\vspace{\abstand}

\begin{bem}
  Let $F=(F_n,u_n)$ be an object of $\pr{\cA}$.
  We can consider $F$ as a map $f:\N_0\rightarrow\Ob(\cA)\times\MorC{\cA}$
  with the property
  \[
    \forall n\in\N_0:\bigl[f(n+1)=\langle A,u\rangle\bigr]\wedge\bigl[f(n)=\langle B,v\rangle\bigr]
    \Rightarrow
    \bigl[v\in\Mor{\cA}{A}{B}\bigr].
  \]
  Then by transfer we see that $\str{F}$ is a projective system
  $(\str{F}_h,\str{u}_h)_{h\in\str{\N_0}}$ in $\str{\cA}$
  with index set $\str{\N_0}$.
  In particular, we have $\str{F}_n=\str{(F_n)}$ for $n\in\N_0$. \\[1 mm]
  For $r\in\str{\N_0}$ we put $\str{F}[h]:=(\str{F}_{h+r})_{h\in\str{\N_0}}$, and by transfer for every $r\geq 0$
  we get a canonical morphism $\str{F}[r]\rightarrow\str{F}$. If $F$ is a zero system with $(F[r]\rightarrow F)=0$
  for a (finite) integer $r\geq 0$, then we also have $(\str{F}[r]\rightarrow\str{F})=0$.
\end{bem}

\vspace{\abstand}

\begin{lemma}\label{lemmaarltrim}
  Let $m\in\str{\N_0}$,
  let $d\in\str{\N}$ be an \emph{infinite} natural number,
  and let $F$ be an Artin-Rees-$l$-adic system in $\pr{\cA}$.
  Then $\image\bigl(\str{F}_{m+d}\rightarrow\str{F}_m\bigr)$ is independent of $d$.
\end{lemma}

\vspace{\abstand}

\begin{proof}
  This follows immediately from \ref{corarml}.
\end{proof}

\vspace{\abstand}

\begin{lemma}\label{lemmaarkertr}
  Let $d_1$ and $d_2$ be \emph{infinite} natural numbers,
  and let $F$ be an Artin-Rees-$l$-adic system in $\pr{\cA}$.
  Then
  \[
    \exists r\in\N_0:
    \image\Bigl(\kernel\bigl[\str{F}_{h+d_1+d_2}\longrightarrow\str{F}_{h+d_1}\bigr]
      \longrightarrow\str{F}_{h+d_1+d_2-r}\Bigr)
    \hookrightarrow
    l^h\cdot\str{F}_{h+d_1+d_2-r}.
  \]
\end{lemma}

\vspace{\abstand}

\begin{proof}
  Because of \ref{corarladic}, there is an $s\in\N_0$ and a short exact sequence
  $0\rightarrow N\rightarrow F[s]\rightarrow G\rightarrow 0$ in $\pr{\cA}$ with a zero system $N$
  and an $l$-adic system $G$.
  Then $N[r]\rightarrow N$ is zero for an $r\geq 0$,
  and \ref{lemmaarker} implies
  \[
    \forall m,n\geq 0:
    \image\Bigl[\kernel\bigl(F_{r+s+m+n}\rightarrow F_{s+n}\bigr)\rightarrow F_{s+m+n}\Bigr]
    \hookrightarrow l^{n+1}\cdot F_{s+m+n}.
  \]
  Transfer of this and specializing to $m:=d_2-r$ and $n:=h+d_1-s$ (note that $m,n\geq 0$, because
  $h$,$d_1$ and $d_2$ are infinite and $r$ and $s$ are finite) gives us
  \begin{multline*}
    \image\Bigl[\kernel\bigl(\str{F}_{r+s+(d_2-r)+(h+d_1-s)}\rightarrow\str{F}_{s+(h+d_1-s)}\bigr)
    \rightarrow\str{F}_{s+(d_2-r)+(h+d_1-s)}\Bigr] \\
    \hookrightarrow l^{(h+d_1-s)+1}\cdot\str{F}_{s+(d_2-r)+(h+d_1-s)},
  \end{multline*}
  and the lemma follows from $h+d_1-s+1\geq h$.
\end{proof}

\vspace{\abstand}

\begin{defi}\label{defiAh}
  Let $N\in\Ns\setminus\N$ be an infinite natural number.
  \begin{enumerate}
    \item
      Because on every object of $\cA$ we have the endomorphism $n$ for every $n\in\N$,
      by transfer we have an endomorphism $N$ for every object of $\str{\cA}$.
    \item
      Let $\str{\cA}_N$ be the full abelian subcategory of $\str{\cA}$ consisting of objects for which
      the endomorphism $N$ is zero.
    \item
      If $A$ is an object of $\str{\cA}$, then we denote the cokernel of $A\xrightarrow{N}A$ by
      $A/N$; obviously $A/N$ is an object of $\str{\cA}_N$.
  \end{enumerate}
\end{defi}

\vspace{\abstand}

\begin{thm}\label{thmarfunctor}
  Let $d_1$ and $d_2$ be two positive infinite numbers in $\str{\N}$.
  For an Artin-Rees-$l$-adic system $F\in\Ob(\pr{\cA})$ put
  \[
    \arfun{h}{F}:=
    \Bigl[\image\bigl(\str{F}_{h+d_1+d_2}\longrightarrow\str{F}_{h+d_1}\bigr)\Bigr]/l^h
    \in\Ob(\str{\cA}_{l^h}).
  \]
  \begin{enumerate}
    \item\label{thmarfunindep}
      $\arfun{h}{F}$ does not depend on the choice of the infinite natural numbers $d_1$ and $d_2$.
    \item\label{thmarfun}
      Let $F\xrightarrow{\tilde{f}}G$ be a morphism in $\artinreesl{\cA}$,
      represented by a morphism $F[r]\xrightarrow{f}G$ in $\pr{\cA}$ for a (finite) $r\geq 0$.
      Then $\str{F}_{h+d_1+d_2}\xrightarrow{\str{f}}\str{F}_{h+d_1+d_2-r}$ induces a morphism
      $\arfun{h}{F}\rightarrow\arfun{h}{G}$ that only depends on $\tilde{f}$
      and which we denote by $\arfun{h}{f}$.
      --- In this way, we get a well-defined covariant additive functor
      $\arfun{h}{}:\artinreesl{\cA}\longrightarrow\str{\cA}_{l^h}$.
    \item\label{thmisoladic}
      If $L$ is an $l$-adic system in $\pr{\cA}$, then
      $\arfun{h}{L}=\str{L}_{h-1}$.
    \item\label{thmarfunex}
      If $F\rightarrow G\rightarrow H\rightarrow 0$ is a sequence in $\artinreesl{\cA}$
      which is exact in $\artinrees{\cA}$,
      then the induced sequence
      \[
        \arfun{h}{F}\rightarrow
        \arfun{h}{G}\rightarrow
        \arfun{h}{H}\rightarrow
        0
      \]
      is exact in $\str{\cA}_{l^h}$ (or -- equivalently -- in $\str{\cA}$).
  \end{enumerate}
\end{thm}

\vspace{\abstand}

\begin{proof}
  To prove \ref{thmarfunindep},
  let $d_1'$, $d_2'$ be another pair of infinite natural numbers. Without loss of generality,
  we can assume $d_1'\leq d_1$, and we note that
  \ref{lemmaarltrim} immediately implies the independence of $d_2$,
  so that we can assume $d_1+d_2=d_1'+d_2'$ as well.
  Then we have an obvious epimorphism
  \[
    \image\bigl(\str{F}_{h+d_1+d_2}\longrightarrow\str{F}_{h+d_1}\bigr)
    \twoheadrightarrow
    \image\bigl(\str{F}_{h+d_1+d_2}\longrightarrow\str{F}_{h+d_1'}\bigr),
  \]
  so we only have to show that modulo $l^h$, this becomes a monomorphism.
  To prove this, we want to assume again that $\str{\cA}$ is an abelian subcategory of $\Abg$. Let
  $x$ be an element of $\str{F}_{h+d_1+d_2}$ whose image in $\str{F}_{h+d_1}$ is mapped to zero in
  $l^h\cdot\image\bigl(\str{F}_{h+d_1+d_2}\longrightarrow\str{F}_{h+d_1'}\bigr)$.
  Then there is $y\in\str{F}_{h+d_1+d_2}$ with
  \[
    x-l^hy\in\kernel\bigl(\str{F}_{h+d_1+d_2}\longrightarrow\str{F}_{h+d_1'}\bigr),
  \]
  and \ref{lemmaarkertr} implies that $x-l^hy$ is mapped to $l^h\cdot\str{F}_{h+d_1}$ (because $d_2$ is infinite,
  and the $r$ whose existence is stated in \ref{lemmaarkertr} is finite) which shows that $x$ is zero in
  $\image\bigl(\str{F}_{h+d_1+d_2}\rightarrow\str{F}_{h+d_1}\bigr)/l^h$.\\[2 mm]
  To prove \ref{thmarfun},
  we have to show that
  there is a morphism $\bar{f}$ that makes the following diagram in $\str{\cA}$ commutative
  (note that $\image(w)=\image(wv)$ because of \ref{lemmaarltrim}):
  \[
    \xymatrix@R=3mm@C=15mm{
      & {\str{F}_{h+d_1+d_2}} \ar[rd]^{f} \ar[dddd]^{u} \ar@{->>}[ldddd]_{u} &
        {\str{G}_{h+d_1+d_2}} \ar[d]_{v} \ar@{->>}[ddddr]^{wv} & \\
      & & {\str{G}_{h+d_1+d_2-r}} \ar[ddd]_{w} \ar@{..>>}[dddr]_{w} & \\
      & & & \\
      & & & \\
      {\image(u)} \ar@{^{(}->}[r] \ar@{->>}[dd]_{p} & {\str{F}_{h+d_1}} & {\str{G}_{h+d_1}} &
        {\image(wv)} \ar@{_{(}->}[l] \ar@{->>}[dd]^{q} \\
      & & & \\
      {\arfun{h}{F}} \ar@{..>}[rrr]_{\bar{f}} & & & {\arfun{h}{G}} \\
    }
  \]
  For simplicity, we want to assume again that $\str{\cA}$ is an abelian subcategory
  of the category of modules over a ring.
  Let $x$ be an element of $\str{F}_{h+d_1+d_2}$ that is mapped to zero in $\arfun{h}{F}$.
  We have to show that $qwf(x)=0\in\arfun{h}{G}$,
  i.e. $wf(x)\in l^h\cdot\image(w)$.\\[1 mm]
  Because of $pu(x)=0$, there is $y\in\str{F}_{h+d_1+d_2}$ with $u(x-l^hy)=0$.
  Then \ref{lemmaarkertr} tells us that there is a $s\in\N_0$ such that
  \[
    \image\bigl[\kernel(u)\longrightarrow\str{F}_{h+d_1+d_2-s}\bigr]
    \hookrightarrow
    l^h\cdot\str{F}_{h+d_1+d_2-s},
  \]
  i.e. if $t$ denotes the map $\str{F}_{h+d_1+d_2}\rightarrow\str{F}_{h+d_1+d_2-s}$,
  there is a $z\in\str{F}_{h+d_1+d_2-s}$ such that $t(x-l^hy)=l^hz$.
  The following diagram is commutative in $\str{\cA}$:
  \[
    \xymatrix{
      {\str{F}_{h+d_1+d_2}} \ar[r]^{f} \ar[d]_{t} & {\str{G}_{h+d_1+d_2-r}} \ar[d] \ar@/^2cm/[dd]^{w} \\
      {\str{F}_{h+d_1+d_2-s}} \ar[r]_{f} & {\str{G}_{h+d_1+d_2-s-r}} \ar[d]^{a} \\
      & {\str{G}_{h+d_1.}} \\
    }
  \]
  We therefore get
  \[
    wf(x)=aft(x)=af\Bigl(l^h\cdot\bigl[z+t(y)\bigr]\Bigr)
    \in l^h\cdot\image(a)
    \stackrel{\ref{lemmaarltrim}}{=}
    l^h\cdot\image(w),
  \]
  which shows that $\bar{f}$ is well-defined.\\[1 mm]
  To see that $\bar{f}$ only depends on $\tilde{f}$,
  we have to show that for every $s\in\N_0$, the morphism $F[r+s]\rightarrow F[r]\xrightarrow{f}G$
  induces the same morphism $\arfun{h}{F}\rightarrow\arfun{h}{G}$ as does $f$.
  But this is clear, because the diagram
  \[
    \xymatrix{
      {\str{F}_{h+d_1+d_2}} \ar[r]^{f} \ar[d] & {\str{G}_{h+d_1+d_2-r}} \ar[d]^{a} \\
      {\str{F}_{h+d_1+d_2-s}} \ar[r]^{f} & {\str{G}_{h+d_1+d_2-(r+s)}} \ar[d]^{b} \\
      & {\str{G}_{h+d_1}} \\
    }
  \]
  is commutative in $\str{\cA}$ and because $\image(ba)=\image(a)$ according to \ref{lemmaarltrim}.
  This completes the proof of \ref{thmarfun}, because the construction is obviously functorial in $\tilde{f}$
  and because $\arfun{h}{}$ obviously is compatible with direct sums.\\[2 mm]
  To prove \ref{thmisoladic}, note that it follows by transfer that if $L$ is $l$-adic, then
  $\str{L}_{h+d_1+d_2}\twoheadrightarrow\str{L}_{h+d_1}$ is an epimorphism and that
  $\str{L}_{h+d_1}/l^h=\str{L}_{h-1}$, i.e.
  \[
    \arfun{h}{L}
    =\Bigl[\image(\str{L}_{h+d_1+d_2}\twoheadrightarrow\str{L}_{h+d_1})\Bigr]/l^h
    =\str{L}_{h+d_1}/l^h
    =\str{L}_{h-1}.
  \]
  For \ref{thmarfunex}, let $F\xrightarrow{f}G\rightarrow H\rightarrow 0$ be exact in
  $\artinrees{\cA}$ with Artin-Rees-$l$-adic systems $F$, $G$ and $H$.
  By definition, we can find a commutative diagram
  \[
    \xymatrix{
      {F} \ar[r]^{f} \ar[d]_{\wr} & {G} \ar[r] \ar[d]_{\wr} & {H} \ar[r] \ar@{=}[d] & 0 \\
      {F'} \ar[r]^{f'} & {G'} \ar[r] & {H} \ar[r] & 0 \\
    }
  \]
  in $\artinrees{\cA}$, whose top row is also exact,
  with $l$-adic systems $F'$ and $G'$.
  Therefore, without loss of generality, we may assume that $F$ and $G$ are already $l$-adic systems.
  Because of \ref{satzprojar}\ref{lemmamorphs},
  we may further assume that $f:F\rightarrow G$ is a morphisms in $\pr{\cA}$.
  Then we may assume that $H$ is the cokernel of $f$,
  and because the cokernel of a morphism of $l$-adic systems is obviously $l$-adic,
  we see that we may assume that
  $F\xrightarrow{f}G\rightarrow H\rightarrow 0$ is an exact sequence of $l$-adic systems in $\pr{\cA}$.
  Now we apply \ref{thmisoladic} which shows that the sequence induced by $\arfun{h}{}$ is simply
  \[
    \str{F}_{h-1}\xrightarrow{(\arfun{h}{f})=f_{h-1}}\str{G}_{h-1}\rightarrow\str{H}_{h-1}\rightarrow 0,
  \]
  which is obviously exact by transfer.
\end{proof}

\vspace{\abstand}

\begin{lemma}\label{lemmaarfact}
  Let $F$ be an Artin-Rees-$l$-adic system in $\pr{\cA}$.
  Then there is an $r\in\N_0$ such that for all $m\in\N_0$,
  the morphism $F_m\rightarrow F_{m-r}$ factorizes over $F_m/l^{m+1}$.
\end{lemma}

\vspace{\abstand}

\begin{proof}
  According to \ref{corarladic}, there is a $s\in\N_0$ and an exact sequence
  $0\rightarrow N\rightarrow F[s]\rightarrow G\rightarrow 0$ in $\pr{\cA}$ with a zero system $N$ and a
  $l$-adic system $G$.
  Let $N[t]\rightarrow N$ be zero.
  For every $n\geq 0$, we then have the following commutative diagram with exact rows in $\cA$:
  \[
    \xymatrix{
      0 \ar[r] & {N_{n}} \ar[r] \ar[d]^{l^{n+1}} & {F_{n+s}} \ar[r] \ar[d]^{l^{n+1}} \ar@{..>}[dl] &
        {G_{n}} \ar[r] \ar[d]^{0} & 0 \\
      0 \ar[r] & {N_{n}} \ar[r] \ar[d]^{0} & {F_{n+s}} \ar[r] \ar[d] & {G_{n}} \ar[r] \ar[d] & {0} \\
      0 \ar[r] & {N_{n-t}} \ar[r] & {F_{n+s-t}} \ar[r] & {G_{n-t}} \ar[r] & {0.} \\
    }
  \]
  The diagram implies that $F_{n+s}\xrightarrow{l^{n+1}}F_{n+s}$ factorizes over $N_n$
  and thus shows that the composition $F_{n+s}\xrightarrow{l^{n+1}}F_{n+s}\rightarrow F_{n+s-t}$ is zero.
  Now chose a $k\geq 0$ satisfying $-k\leq s-t\leq s\leq k$, then we have the following commutative diagram in $\cA$:
  \[
    \xymatrix{
      {F_{n+k}} \ar[r]^{l^k} \ar@{..>}@/_10mm/[dddrr]_{0} \ar@{..>}@/^8mm/[rr]^{l^{(n+k)+1}} &
        {F_{n+k}} \ar[d] \ar[r]^{l^{n+1}} & {F_{n+k}} \ar[d] \ar@{..>}[dddr] & \\
      & {F_{n+s}} \ar[rd]_{0} \ar[r]^{l^{n+1}} & {F_{n+s}} \ar[d] & \\
      & & {F_{n+s-t}} \ar[d] & \\
      & & {F_{n-k}} \ar@{=}[r] & {F_{(n+k)-2k}}\\
    }
  \]
  So if we choose $r:=2k$ and substitute $m-k$ for $n$, the claim follows.
\end{proof}

\vspace{\abstand}

\begin{thm}\label{thmfaithful}
  \mbox{}\\[-5 mm]
  \begin{enumerate}
    \item\label{satzstrfun}
      We get an exact functor
      $\artinrees{\cA}\rightarrow\artinrees{\str{\cA}}$
      by sending a system $F$ to the system $(\str{F}_n)_{n\in\N_0}$,
      and this functor maps Artin-Rees-$l$-adic systems to Artin-Rees-$l$-adic systems and thus
      induces a functor $\arstrfun{}:\artinreesl{\cA}\rightarrow\artinreesl{\str{\cA}}$.
    \item\label{satzzlhfun}
      We get a right exact functor
      $\arzlhfun{h}{}:\str{\cA}_{l^h}\rightarrow\artinreesl{\str{\cA}}$ by defining
      $\arzlhfun{h}{A}:=(A/l^{n+1})_{n\in\N_0}$ for $A\in\Ob(\str{\cA}_{l^h})$.
    \item\label{satzisofun}
      There is a canonical isomorphism
      $\varphi:\arzlhfun{h}{}\circ{}\arfun{h}{}\xrightarrow{\sim}\arstrfun{}$
      of functors from $\artinreesl{\cA}$ to $\artinreesl{\str{\cA}}$.
  \end{enumerate}
\end{thm}

\vspace{\abstand}

\begin{proof}
  \ref{satzstrfun} is immediately clear by transfer,
  and \ref{satzzlhfun} is trivial.
  To prove \ref{satzisofun},
  let $d_1$ and $d_2$ be two infinite natural numbers,
  and let $F$ be an arbitrary Artin-Rees-$l$-adic system in $\pr{\cA}$.
  We first want to show that for every $n\in\N_0$,
  the morphism $\str{F}_{h+d_1+d_2}\rightarrow\str{F}_n$ induces a morphism
  $\psi:\arfun{h}{F}/l^{n+1}\rightarrow\str{F}_n/l^{n+1}$ which is functorial in $F$.
  Denote the morphism $\str{F}_{h+d_1+d_2}\rightarrow\str{F}_{h+d_1}$ by $u$.
  Because of \ref{lemmaarkertr}, we have the following commutative diagram in $\str{\cA}$ with exact rows which
  induces a morphism $\vartheta:\image(u)\rightarrow\str{F}_n/l^h$:
  \[
    \xymatrix{
      0 \ar[r] & {\kernel(u)} \ar[r] \ar[d] & {\str{F}_{h+d_1+d_2}} \ar[r]^{u} \ar[d] &
        {\image(u)} \ar[r] \ar@{..>}[d]^{\vartheta} & 0 \\
      0 \ar[r] & {l^h\cdot\str{F}_n} \ar[r] & {\str{F}_n} \ar[r] & {\str{F}_n/l^h} \ar[r] & {0,} \\
    }
  \]
  and then the following commutative diagram with exact rows in $\str{\cA}$ induces $\psi$,
  and the construction is obviously functorial in $F$:
  \[
    \xymatrix{
      {\image(u)} \ar[r]^{l^{n+1}} \ar[d] & {\image(u)} \ar[r] \ar[d] &
        {\arfun{h}{F}/l^{n+1}} \ar[r] \ar@{..>}[d]^{\psi} & 0 \\
      {\str{F}_n/l^h} \ar[r]^{l^{n+1}} & {\str{F}_n/l^h} \ar[r] & {\str{F}_n/l^{n+1}} \ar[r] & {0.} \\
    }
  \]
  Combining this result with \ref{lemmaarfact}, we see that there is an $r\in\N_0$ such that for every $n\in\N_0$,
  we have a morphism $\varphi_{F,n}:\arfun{h}{F}/l^{n+1}\rightarrow\str{F}_{n-r}$ in $\str{\cA}$
  which induces a morphism $\bigl(\arfun{h}{F}/l^{n+1}\bigr)[r]\rightarrow\bigl(\str{F}_n\bigr)$
  in $\pr{\str{\cA}}$ and thus a morphism $\varphi_F:\arzlhfun{h}{\arfun{h}{F}}\rightarrow\arstrfun{F}$
  in $\artinreesl{\str{\cA}}$ which is obviously functorial in $F$. The system of all $(\varphi_F)$ therefore
  defines a morphism of functors $\arzlhfun{h}{}\circ\arfun{h}{}\xrightarrow{\varphi}\arstrfun{}$.\\[2 mm]
  To see that $\varphi$ is actually an \emph{isomorphism} of functors, we have to show that $\varphi_F$
  is an isomorphism in $\str{\cA}$ for every Artin-Rees-$l$-adic system $F$.
  By definition, there is an $l$-adic system $G$ in $\pr{\cA}$ and an isomorphism $f:F\xrightarrow{\sim}G$
  in $\artinrees{\cA}$.
  Then because $\varphi$ is a morphism of functors, the following diagram in $\str{\cA}$ is commutative:
  \[
    \xymatrix{
      {\arzlhfun{h}{\arfun{h}{F}}} \ar[d]^{\wr}_{\arzlhfun{h}{\arfun{h}{f}}} \ar[r]^{\varphi_F} &
        {\arstrfun{F}} \ar[d]_{\wr}^{\arstrfun{f}} \\
      {\arzlhfun{h}{\arfun{h}{G}}} \ar[r]_{\varphi_G} & {\arstrfun{G}} \\
    }
  \]
  This shows that without loss of generality, we can assume that $F$ is an $l$-adic system.
  But then
  \[
    \arzlhfun{h}{\arfun{h}{F}}
    \stackrel{\ref{thmarfunctor}\ref{thmisoladic}}{=}
    \arzlhfun{h}{\str{F}_{h-1}}
    =\bigl(\str{F}_{h-1}/l^{n+1}\bigr)
    =\bigl(\str{F}_n\bigr)
    =\arstrfun{F}.
  \]
\end{proof}

\vspace{\abstand}

\begin{cor}\label{corfaithful}
  The functors $\arfun{h}{}:\artinreesl{\cA}\rightarrow\str{\cA}_{l^h}$
  and $\arstrfun{}:\artinreesl{\cA}\rightarrow\artinreesl{\str{\cA}}$
  are faithful and reflect isomorphisms.
\end{cor}

\vspace{\abstand}

\begin{proof}
  We only have to prove this for $\arstrfun{}$,
  because then \ref{thmfaithful} immediately implies it for $\arfun{h}{}$ as well.
  It is very easy to see that $\arstrfun{}$ is faithful, so let $f:F\rightarrow G$ be a morphism
  in $\artinreesl{\cA}$ such that $\arstrfun{f}$ is an isomorphism.
  By definition, there is an $r\in\N_0$ and an exact sequence
  $0\rightarrow M\rightarrow F[r]\xrightarrow{\bar{f}}G\rightarrow N\rightarrow 0$
  in $\pr{\cA}$ such that $\arstrfun{M}$ and $\arstrfun{N}$ are zero systems in $\pr{\str{\cA}}$.
  But transfer then immediately implies that $M$ and $N$ are also zero systems, i.e. $f$ is an isomorphism.
\end{proof}

\vspace{\abstand}

\begin{lemma}\label{lemmansprojlim}
  In the special case where $\cA$ is equivalent to the category $\fin{\Abe}$ of finite abelian groups,
  let $(A_n)_{n\in\N_0}$ be a projective system of finite abelian groups, considered as objects of $\cA$.
  Then we have a canonical epimorphism
  $A_h\twoheadrightarrow\varprojlim_{n\in\N_0}A_n$ in $\str{\cA}\subseteq\Abe$
  (we consider $\str{A}$ as an exact subcategory of $\Abe$ because of \ECref{bemrmodfin}),
  and an $a\in A_h$ is in the kernel
  if and only if there is an infinite $h'\leq h$ such that $a$ maps to zero under the canonical morphism
  $A_h\rightarrow A_{h'}$.
\end{lemma}

\vspace{\abstand}

\begin{proof}
  We consider the projective system as a covariant functor $F$ from the cofiltered category $\N_0$ into the
  category $\cA$,
  so that we get a functor $\str{F}:\str{\N_0}\rightarrow\str{\cA}\subseteq\Abe$.
  In particular, we get $A_h:=(\str{F})(h)$ and for all $n\leq h$ compatible maps $A_h\rightarrow(\str{F})(n)$.
  If $n$ is finite, we have $(\str{F})(n)=A_n$ because $A_n$ is \emph{finite}
  and because of \ECrefs{satzrmodfun}{satzrmodfuni},
  so that we get a canonical map
  $A_h\rightarrow\varprojlim_{n\in\N_0}A_n$ as claimed.\\[2 mm]
  To see that this map is surjective, let $a=(a_n)$ be an element of $\varprojlim_{n\in\N_0}A_n$.
  If $\cB$ is an $\hat{S}$-small full abelian subcategory of $\Abe$ that contains $\cA$ as a full subcategory
  and the group $\Z$ as an object,
  then we can consider $a$ as a map $\N_+\rightarrow\MorC{\cD}$ with the property
  \[
    \Bigl(\forall n\in\N_+:a(n)\in\Mor{\cD}{\Z}{A_n}\Bigr)
    \wedge
    \Bigl(\forall n,n'\in\N_+:\bigl[n>n'\Rightarrow a(n')=F(n\rightarrow n')\circ a(n)\bigr]\Bigr),
  \]
  and transfer of this shows that $\str{a}$ is an element of
  $\varprojlim_{n\in\str{\N_+}}\Mor{\str{\cB}}{\str{\Z}}{A_n}
  \stackrel{\text{\tiny\ECrefs{satzrmodfun}{satzrmodfuni}}}{=}
  \varprojlim_{n\in\str{\N_+}}A_n$, i.e
  $a_h:=(\str{a})(h)$ is a preimage of $a$.\\[2 mm]
  To prove the statement about the kernel, let $b$ be an element of $A_h$,
  i.e. a morphism from $\str{\Z}$ to $A_h$ in $\str{\cB}$ according to \ECrefs{satzrmodfun}{satzrmodfuni}.
  If $b$ maps to zero in $A_{h'}$ for an infinite $h'\leq h$, then it also maps to zero in all $A_n$ for $n\in\N_+$
  (because those are all smaller than $h'$),
  so that $b$ maps to zero in $\varprojlim_{n\in\N_+}A_n$.\\[1 mm]
  Now let $b$ be in the kernel of $A_h\rightarrow\varprojlim_{n\in\N_+}A_n$.
  Then the set
  \[
    N:=\Bigl\{n\in\str{\N_+}^{\leq h}\Bigr\vert\str{\Z}\xrightarrow{b}A_h\rightarrow A_n=0\Bigr\}
  \]
  is an \emph{internal} subset\footnote[2]{Note that $b$ and $A_h\rightarrow A_n$ are internal because
  they are elements of $\MorC{\str{\cB}}=\str{\MorC{\cB}}$!}
  of $\str{\N_+}^{\leq h}$ that contains the \emph{external} set $\N_+$ and therefore
  must contain an element $h'$ that is not in $\N_+$ and that hence is infinite. This completes the proof.
\end{proof}

\vspace{\abstand}

\begin{cor}\label{corzlmod}
  Let $*:\hat{S}\rightarrow\widehat{\str{S}}$ be an enlargement,
  let $\cA$ be an $\hat{S}$-small abelian category that is equivalent to the category
  of finite abelian groups,
  let $l$ be a prime number,
  let $h$ be an infinite natural number,
  and let $F$ be an Artin-Rees-$l$-adic system in $\cA$.
  \begin{enumerate}
    \item\label{corzlhfin}
      The system $\arzlhfun{h}{\arfun{h}{F}}$ is an $l$-adic system in $\pr{\cA}$ which is
      canonically isomorphic to $F$ in $\artinrees{\cA}$.
    \item\label{corzliso}
      Let $d_1$ and $d_2$ be infinite natural numbers. Then
      we have canonical isomorphisms of finitely generated $\Z_l$-modules
      \[
        \varprojlim_{n\in\N_0}F_n
        =\arfun{h}{F}\otimes_{\str{\Z/l^h}}\Z_l
        =\image(\str{F}_{d_1+d_2}\rightarrow\str{F}_{d_1})\otimes_{\str{\Z}}\Z_l.
      \]
      (Note that we know from \ECref{satzrmodfun} that the objects of $\str{\cA}$ are $\str{\Z}$-modules,
      and that $\Z_l$ is a quotient of $\str{\Z}$ according to \ECrefs{bsplim}{bspprojlim}
  \end{enumerate}
\end{cor}

\vspace{\abstand}

\begin{proof}
  $\arzlhfun{h}{\arfun{h}{F}}$ is obviously an $l$-adic system in $\pr{\str{\cA}}$,
  and according to \ref{thmfaithful}\ref{satzisofun}, it is canonically isomorphic to $\arstrfun{F}$.
  Because of \ref{corarladic}, there then is an $r\in\N_0$ and an exact sequence
  \[
    0\rightarrow N\rightarrow\arstrfun{F}[r]\rightarrow\arzlhfun{h}{\arfun{h}{F}}\rightarrow 0
  \]
  in $\pr{\str{\cA}}$ with a zero system $N$.
  Because all the $\str{F}_n=F_n$ are finite abelian groups,
  we see that all $N_n$ and $\arzlhfun{h}{\arfun{h}{F}}_n$ must also be finite
  (compare \ECrefs{satzrmodfun}{satzrmodfuniii}),
  and this implies that
  the same sequence can also be considered as an exact sequence in $\pr{\cA}$. This proves \ref{corzlhfin}.\\[2 mm]
  It is well known that $(A_n)\mapsto\varprojlim_{n\in\N_+}A_n$ establishes an equivalence between $\artinreesl{\cA}$
  and the category of finitely generated $\Z_l$-modules (compare \cite[\S12]{kiehl1}).
  Because of \ref{corzlhfin} we therefore only have to prove that the kernel of the
  epimorphism $\arfun{h}{F}\twoheadrightarrow\varprojlim_{n\in\N_+}\arfun{h}{F}/l^n$
  from \ref{lemmansprojlim} equals $l^{\leq h}\cdot\arfun{h}{F}$ where $l^{\leq h}$ by definition is
  the $\str{\Z}$-ideal
  $\langle l^{h'}\vert h'\in\str{\N_+}^{\leq h}\setminus\N_+\rangle$.
  But this follows immediately from the description of the kernel in \ref{lemmansprojlim}.
\end{proof}

\vspace{\abstand}

\subsection{Application to $l$-adic cohomology}

\mbox{}\\
\vspace{\abstand}

As in the first chapters,
let $S$ be a noetherian scheme,
let $U$ be a universe such that the category $\text{PreShv}(Et(X), U-Sets)$ contains all representable presheaves
for all $X\in Sch/S$,
let $\hat{S}$ be a superstructure that contains all occurring categories,
and let $*:\hat{S}\rightarrow\widehat{\str{S}}$ be an enlargement.\\

Let $X$ be an object of $Sch/S$, and let $l\in\N$ be a prime number.
Then by definition, an \emph{$l$-adic sheaf on $X$} is an $l$-adic system in the category
$\artinrees{\constr{X}}$,
and an \emph{Artin-Rees $l$-adic sheaf on $X$} (or \emph{AR-$l$-adic sheaf} for short)
is an object of the category $\artinreesl{\constr{X}}$.
As already mentioned above, the category $\artinreesl{\constr{X}}$ of all AR-$l$-adic sheaves on $X$
is a full exact abelian subcategory of $\artinrees{\constr{X}}$
--- for all this compare \cite[\S12]{kiehl1}.\\

Let $h\in\str{\N}\setminus\N$ be an infinite natural number.
Then obviously the category $\str{\constr{X}}_{l^h}$ from \ref{defiAh}
equals the category $\str{\modconstr{X}{\str{\Z}/l^h}}$ defined in chapter \ref{chone}, so that
we get the following corollary from \ref{thmarfunctor} and \ref{corfaithful}:

\vspace{\abstand}

\begin{cor}\label{corladic}
  Let $X$ be an object of $Sch/S$,
  let $l\in\N$ be a prime number,
  and let $h\in\str{\N}\setminus\N$ be an infinite natural number.
  Then
  \begin{enumerate}
    \item\label{corladici}
      $\arfun{h}$ is a faithful, right exact functor from the category of AR-$l$-adic sheaves on $X$ to the
      category $\str{\modconstr{X}{\str{\Z}/l^h}}$ which reflects isomorphisms.
    \item\label{corladicii}
      If $\cF=(\cF_n)$ is an $l$-adic sheaf on $X$,
      then $\arfun{h}\cF=\str{\cF}_{h-1}$.
  \end{enumerate}
\end{cor}

\vspace{\abstand}

\noindent
Now consider the special case where $S=\spec{k}$ is the spectrum of a field $k$.
Then we get:

\vspace{\abstand}

\begin{thm}\label{thmladic}
  Let $X$ be a variety over $k$,
  and let $i\in\N_0$ be a natural number.
  \begin{enumerate}
    \item\label{thmladici}
      For each prime number $l\in\N$, and for two infinite natural numbers $d_1$ and $d_2$,
      we have canonical isomorphisms
      \[
        \begin{array}{cccc}
          \image\bigl[H^i_c(\bar{X},\str{\Z}/l^{d_1+d_2})\rightarrow H^i_c(\bar{X},\str{\Z}/l^{d_1})\bigr]
            \otimes_{\str{\Z}}\Z_l &
          \xrightarrow{\sim} &
          H^i_c(\bar{X},\Z_l) &
          \text{and} \\[3 mm]
          \image\bigl[H^i_c(\bar{X},\str{\Z}/l^{d_1+d_2})\rightarrow H^i_c(\bar{X},\str{\Z}/l^{d_1})\bigr]
            \otimes_{\str{\Z}}\Q_l &
          \xrightarrow{\sim} &
          H^i_c(\bar{X},\Q_l). &
          \\
        \end{array}
      \]
    \item\label{thmladicii}
      If $X$ is smooth and projective, then
      for \emph{almost all} prime numbers $l\in\N$, we have canonical isomorphisms
      \[
        \begin{array}{cccc}
          H^i(\bar{X},\str{\Z}/l^h)\otimes_{\str{\Z}}\Z_l &
          \xrightarrow{\sim} &
          H^i(\bar{X},\Z_l) &
          \text{and} \\[3 mm]
          H^i(\bar{X},\str{\Z}/l^h)\otimes_{\str{\Z}}\Q_l &
          \xrightarrow{\sim} &
          H^i(\bar{X},\Q_l). &
          \\
        \end{array}
      \]
  \end{enumerate}
\end{thm}

\vspace{\abstand}

\begin{proof}
  Let $k_s$ be a separable closure of $k$,
  let $f:\bar{X}\rightarrow\spec{k_s}$ be the structure morphism of $\bar{X}=X\otimes_kk_s$,
  and let $\cF$ be the $l$-adic sheaf $(\Z/l^{n+1})_n$ on $\bar{X}$.
  Because of the finiteness theorem for $l$-adic sheaves (\cite[12.15]{kiehl1}),
  the system $\cG:=(\RR^if_!\cF_n)_n=(H^i_c(\bar{X},\Z/l^{n+1})_n$ is then an AR-$l$-adic sheaf on $\spec{k}$.
  But a constructible sheaf on $\spec{k}$ is just a finite abelian group with an action of the absolute Galois group $G_k$,
  so \ref{thmladici} follows from \ref{corzlmod}\ref{corzliso}.\\[1 mm]
  To prove \ref{thmladicii},
  because of \ref{thmladici} and \ref{corladic}\ref{corladicii}
  it suffices to show that $\cG$ is not only an AR-$l$-adic sheaf, but an $l$-adic sheaf.
  For that we use Gabber's already cited result from \cite{gabber}
  that for almost all primes $l$, the $l$-adic cohomology of a smooth projective variety has no torsion.
  Let $l$ be one of those $l$ for which $H^{i+1}(\bar{X},\Z_l)$ is torsion free.
  We claim that in this case $\cG$ is an $l$-adic sheaf.
  We only have to show that for all $n\in\N_+$, we have $H^i(\bar{X},\Z/l^{n+1})/l^n\cong H^i(\bar{X},\Z/l^n)$,
  because the other conditions from \ref{defiar}\ref{defiladic}
  are trivially satisfied.
  So let $n\in\N_+$ be a positive natural number.
  As in \cite[V.1.11]{milne},
  we have for all $m\in\N_+$ the following short exact sequence of finite abelian groups:
  \[
    0\longrightarrow
    H^i(\bar{X},\Z_l)/l^m\longrightarrow
    H^i(\bar{X},\Z/l^m)\longrightarrow
    H^{i+1}(\bar{X},\Z_l)_{l^m}\longrightarrow
    0,
  \]
  but by hypothesis there is no $l^m$-torsion in $H^{i+1}(\bar{X},\Z_l)$,
  so we get an isomorphism $H^i(\bar{X},\Z_l)/l^m\cong H^i(\bar{X},\Z/l^m)$ for all $m\in\N_+$.
  But then
  \[
    H^i(\bar{X},\Z/l^{n+1})/l^n
    \cong\bigl[H^i(\bar{X},\Z_l)/l^{n+1}\bigr]/l^n
    \cong H^i(\bar{X},\Z_l)/l^n
    \cong H^i(\bar{X},\Z/l^n).
  \]
  This completes the proof of the theorem.
\end{proof}

\vspace{\abstand}

\begin{appendix}
  \section{Enlargements of superstructures and categories}

\mbox{}\\
\vspace{\abstand}

In this chapter we want to recall and  summarize the most important
definitions from the theory of enlargements of superstructures
and categories. For details we refer to \cite{loeb}, \cite{lanrog} and \cite{enlcat}.

\vspace{\abstand}

For a set $M$, let $\cP(M)$ denote the power set of $M$, i.e. the set of all subsets of $M$.

\vspace{\abstand}

\begin{defi}(Superstructure)\\
  Let $S$ be an infinite set whose elements are no sets.
  Such a set we call a \emph{base set}, and its elements we call \emph{base elements}.
  We define $\hat{S}$, the \emph{superstructure over $S$}, as follows:
  \[
    \hat{S}:=\bigcup_{n=0}^\infty S_n
    \;\;\;\;\;\;
    \text{where $S_0:=S$ and $\forall n\geq 1:S_n:=S_{n-1}\cup\cP(S_{n-1})$.}
  \]
\end{defi}

\vspace{\abstand}

In the superstructure $\hat{S}$ to a base set $S$, we will find most of the mathematical objects of interest
related to $S$:
First of all, for sets $M,N\in\hat{S}$, the product set $M\times N$ is again an element of $\hat{S}$ when
we identify an ordered pair $\langle a,b\rangle$ for $a\in M$, $b\in N$ with the set
$\{a,\{a,b\}\}$,
and for sets $M_1,\ldots,M_n\in\hat{S}$, the product set
$M_1\times\ldots\times M_n:=(M_1\times\ldots\times M_{n-1})\times M_n$
is an element of $\hat{S}$.
Therefore, relations between two sets $M,N\in\hat{S}$ and in particular functions from $M$ to $N$
are again elements of $\hat{S}$.

For example, if $S$ contains the set of real numbers $\R$, then $\hat{S}$ will contain the sets $\R^n$ for
$n\in\N_+$ as well as functions between subsets of $\R^n$ and $\R^m$, the set of continuous functions between such sets
or the set of differentiable functions and so on.

\vspace{\abstand}

\begin{defi}\label{defiSsmallcat}
  A small category $\cC$ is called $\hat{S}$--small if the set $\MorC{\cC}:=\bigsqcup_{X,Y\in\Ob(\cC)}\Mor{\cC}{X}{Y}$
  is an Element in $\hat{S}$.
\end{defi}

\vspace{\abstand}

\begin{defi}\label{defienlargement}(Enlargement) \\
  Let $*:\hat{S}\rightarrow\hat{W}$ be a map between superstructures.
  For $\tau\in\hat{S}$ we denote the image of $\tau$ under $*$ by $\str{\tau}$,
  and for a formula $\varphi$ in $\hat{S}$,
  we define $\str{\varphi}$ to be the formula in $\hat{W}$ that we get when we replace any constant
  $\tau$ occurring in $\varphi$ by $\str{\tau}$.\\[1mm]
  We call $*$ an \emph{enlargement} if the following conditions hold:
  \begin{enumerate}
    \item
      $\str{S}=W$.
      (Because of this property, we will often write $*:\hat{S}\rightarrow\widehat{\str{S}}$.)
    \item(\emph{transfer principle})\\
      If $\varphi$ is a statement in $\hat{S}$, then $\varphi$ is true
      iff $\str{\varphi}$ is true.
      (A \emph{statement} in a superstructure is a statement build from terms --- using only elements of the superstructure
      as constants --- and the usual logical connectives, with universal and existential quantification over
      sets which are elements of the superstructure.)
    \item(\emph{saturation principle})\label{saturation}\\
      Put $\cI:=\bigcup_{A\in\hat{S}\setminus S}\str{A}\subseteq\hat{W}$,
      let $I$ be a nonempty set whose cardinality is not bigger than that of $\hat{S}$,
      and let $\{U_i\}_{i\in I}$ be a family of nonempty sets $U_i\in\cI$
      with the property that for all \emph{finite} subsets $ J\subseteq I$,
      the intersection $\bigcap_{j\in J}U_j$ is nonempty. Then $\bigcap_{i\in I}U_i\neq\emptyset$.
  \end{enumerate}
  If $*$ is an enlargement, we call the elements of $\cI$ the \emph{internal elements} of $\hat{W}$.
\end{defi}

\vspace{\abstand}

\begin{thm}
  For any base set $S$, there exists an enlargement $*:\hat{S}\rightarrow\widehat{\str{S}}$.
\end{thm}

\vspace{\abstand}

\begin{satzdefi}(Enlargements of categories and functors, compare \ECref{satzstarcat} and \ECref{satzfunctors})\\
  Let $*:\hat{S}\ra\widehat{\str{S}}$ be an enlargement,
  and let $\cC$ be an $\hat{S}$-small category.
  Then the enlargement $\str{\MorC{\cC}}$ of the set of morphisms in $\cC$
  is the set of morphisms of an $\widehat{\str{S}}$-small category $\str{\cC}$ in a natural way,
  and we call $\str{\cC}$ the \emph{enlargement of $\cC$}.\\[1 mm]
  If ${\cD}$ is another $\hat{S}$-small category,
  and if $F:\cC\ra\cD$ is a functor,
  then the enlargement $\str{F}$ of $F$ (where we consider $F$ as a map from $\MorC{\cC}$ to $\MorC{\cD}$)
  is a functor from $\str{\cC}$ to $\str{\cD}$
  which we call the \emph{enlargement of $F$}.
\end{satzdefi}

\vspace{\abstand}

\begin{satzdefi}
  Let $\hat{S}$ be a superstructure.
  An additive, abelian or triangulated category is called \emph{$\hat{S}$-small}
  if the underlying category is $\hat{S}$-small;
  an (additive, abelian, triangulated,...) fibration $\cC\ra\cD$ is called \emph{$\hat{S}$-small}
  if $\cC$ and $\cD$ are $\hat{S}$-small.\\[1 mm]
  If $*:\hat{S}\ra\widehat{\str{S}}$ be an enlargement,
  then the enlargement of an $\hat{S}$-small additive, abelian or triangulated category
  is an $\widehat{\str{S}}$-small additive, abelian respectively triangulated category,
  the enlargement of an $\hat{S}$-small (additive, abelian, triangulated,...) fibration
  is an $\widehat{S}$-small (additive, abelian, triangulated,...) fibration,
  and the enlargement of an additive, left exact, right exact or exact functor
  is again additive, left exact, right exact respectively exact.
\end{satzdefi}

\vspace{\abstand}

\end{appendix}

\bibliographystyle{alpha}
\bibliography{Literatur}

\end{document}